\documentclass[12pt,leqno,oneside]{amsart}

\usepackage{mathrsfs}
\usepackage{amsmath,amstext,amsthm,amssymb}
\usepackage{charter}
\usepackage{mathtools}
\usepackage{microtype}
\usepackage{enumitem}
\usepackage[a4paper,left=2.35cm,right=2.35cm,top=2.25cm,bottom=2.35cm,footskip=1cm]{geometry}
\usepackage[hidelinks]{hyperref}
\hypersetup{
  pdftitle={Mass equidistribution and fluctuations for random holomorphic sections of big line bundles with continuous metrics},
  pdfauthor={Turgay Bayraktar and Afrim Bojnik}
}

\numberwithin{equation}{section}

\newtheorem{Theorem}{Theorem}[section]
\newtheorem{Proposition}[Theorem]{Proposition}
\newtheorem{Lemma}[Theorem]{Lemma}
\newtheorem{Corollary}[Theorem]{Corollary}
\theoremstyle{definition}

\newtheorem{Example}[Theorem]{Example}
\theoremstyle{remark}
\newtheorem{Remark}[Theorem]{Remark}

\newcommand{\imat}{\sqrt{-1}}
\newcommand{\norm}[1]{\lVert #1\rVert}

\newcommand{\ov}{\overline}
\newcommand{\Tr}{\operatorname{Tr}}
\newcommand{\Var}{\operatorname{Var}}
\newcommand{\Cov}{\operatorname{Cov}}
\newcommand{\supp}{\operatorname{supp}}
\newcommand{\vol}{\operatorname{vol}}
\newcommand{\C}{\mathbb{C}}
\newcommand{\R}{\mathbb{R}}

\newcommand{\E}{\mathbb{E}}
\renewcommand{\P}{\mathbb{P}}
\newcommand{\1}{\mathbf{1}}
\newcommand{\cH}{\mathscr{H}}
\newcommand{\cM}{\mathscr{M}}

\begin{document}

\title{Mass equidistribution and fluctuations for random holomorphic sections of big line bundles with continuous metrics}

\author{Turgay Bayraktar}
\address{Faculty of Engineering and Natural Sciences, Sabanc{\i} University, \.{I}stanbul, Turkey}
\email{tbayraktar@sabanciuniv.edu}
\thanks{T. Bayraktar is partially supported by T\"{U}B\.{I}TAK grant ARDEB-1001/124F370.}

\author{Afrim Bojnik}
\address{Faculty of Mathematics and Computer Science,
Jagiellonian University,
\L{}ojasiewicza 6,
30-348 Krak\'ow, Poland}
\email{afrim.bojnik@uj.edu.pl}
\thanks{A. Bojnik is partially supported
by the National Science Centre, Poland, under the
Weave-UNISONO grant no.~2025/07/Y/ST1/00146.}


\date{}

\begin{abstract}
We study mass equidistribution and fluctuations of random
holomorphic sections of big line bundles over compact complex
manifolds, equipped with continuous Hermitian metrics and
Bernstein--Markov measures. For arrays of centered coefficients
that are independent and identically distributed within each row,
we establish sharp moment conditions for mass equidistribution,
covering a broad class of distributions, including Gaussian
and certain discrete laws. On compact K\"ahler manifolds,
we also show that mass equidistribution implies equidistribution
of zeros under a natural condition on the equilibrium measure.
We further establish central limit theorems beyond the Gaussian
setting, show that fluctuations may depend on the choice of
orthonormal basis, and determine their limiting variances
for Haar-almost every sequence of bases.

\end{abstract}

\maketitle

\section{Introduction}\label{S:Intro}

The asymptotic distribution of the mass measures of holomorphic sections is a natural problem in complex geometry. It describes how the mass of these sections is distributed over the underlying manifold. For sections of increasing tensor powers of a line bundle, one asks whether their suitably normalized  mass measures converge to an equilibrium measure determined by the geometry.

In the positive setting, this problem is closely related to
quantum ergodicity and the distribution of zeros of holomorphic
sections; see Shiffman--Zelditch \cite{ShiffmanZelditch1999}
and Zelditch \cite{Zelditch2018}. For arbitrary sequences of sections, however, such convergence cannot be expected without further assumptions. This naturally
leads to the study of random sections, where one seeks to
understand under what probabilistic conditions the geometry
determines the asymptotic behavior of individual sections. However, much of the existing theory concerns Gaussian ensembles, whose invariance under unitary changes of basis
makes them particularly well suited to geometric analysis.
Bayraktar \cite{Bayraktar2020} extended mass equidistribution
to random polynomials with sub-Gaussian coefficients. It is
therefore natural to ask to what extent these phenomena persist
under weaker assumptions on the coefficient distributions,
and whether they continue to hold in more general geometric
settings.

We address these questions for random holomorphic sections
of powers of a big line bundle equipped with a continuous
Hermitian metric. Working with inner products associated with Bernstein--Markov measures, we consider independent coefficients whose distributions need not be
Gaussian or even sub-Gaussian. We characterize mass equidistribution
in probability and identify a sharp moment condition for
almost-sure convergence. Thus, the results encompass the
Gaussian and sub-Gaussian cases while requiring considerably
weaker moment assumptions. We also establish a
deterministic principle relating mass equidistribution to
the distribution of zeros in the K\"ahler case under very natural assumptions.

Beyond equidistribution, we investigate the fluctuations of
the mass measures. A central limit theorem for Gaussian
random sections was established in \cite{BayraktarBojnik2026}
for sequences of positive line bundles, using the classical
Sodin--Tsirelson \cite{SodinTsirelson2004}
theorem for nonlinear functionals of Gaussian processes.
For related results on the distribution and fluctuations
of zeros of random holomorphic sections associated with
sequences of line bundles, see e.g.\
\cite{Bay, BCHM, BayraktarComanMarinescu2020,BojnikGunyuz2025,BojnikGunyuz2026} and the references therein.
Our approach here takes a different route. By expressing mass
statistics as quadratic forms associated with Toeplitz
operators, we connect the geometric problem with random quadratic forms.
This allows us to establish central limit theorems for a broader
class of coefficient distributions, including Gaussian laws as a special case.
For general coefficient laws, however, the limiting
fluctuations may depend on the choice of orthonormal basis.
This leads us to investigate
their behavior for typical bases. By equipping the unitary
groups parametrizing orthonormal bases with Haar probability
measures, we also establish central limit theorems with explicit
limiting variances for Haar-almost every sequence of
orthonormal bases.

\subsection{Random sections and mass measures}

Let $X$ be a compact connected complex manifold of dimension $n\geq1$,
and let $L\to X$ be a big holomorphic line bundle with a continuous
Hermitian metric $h$. Let $\mu$ be a Bernstein--Markov probability
measure for $(K,h)$, where $K=\supp\mu$ is non-pluripolar. We write
$\mu_{\rm eq}$ for the associated weighted equilibrium probability
measure. The relevant geometric notions are recalled in
Section~\ref{S:Geometry}.

For $k\geq1$, set $\cH_k=H^0(X,L^k)$ and $N_k=\dim\cH_k$. Since $L$ is big, as $k\to\infty$,
\begin{equation}\label{e:dimension-asymptotic-intro}
N_k=\frac{\vol(L)}{n!}k^n+o(k^n).
\end{equation}
We henceforth work with sufficiently large $k$, so that $N_k>0$. Equip $\cH_k$
with the following Hermitian inner product
\begin{equation}\label{e:L2-intro}
(s,t)_k:=\int_K h^k(s,t)\,d\mu.
\end{equation}
Throughout, Hermitian inner products are linear in the first variable.
For each $k$, fix an orthonormal basis
$\{S_{1,k},\ldots,S_{N_k,k}\}$ of $\cH_k$.

Let $(\Omega,\mathcal F,\P)$ be a probability space, and let
$\xi$ be a complex random variable satisfying
$\E\xi=0$ and $\P(\xi\ne0)>0$. For each $k$,
let $\xi_{1,k},\ldots,\xi_{N_k,k}$ be independent copies of $\xi$,
all defined on the same probability space. We consider random
holomorphic sections of the form
\begin{equation}\label{e:random-section-intro}
s_k(\omega)=\sum_{j=1}^{N_k}\xi_{j,k}(\omega)S_{j,k},
\qquad \omega\in\Omega.
\end{equation}
Unless otherwise stated, the orthonormal bases are fixed and
deterministic, and no independence between different rows
of coefficients is assumed.

The \emph{dimension-normalized mass measure} and, for nonzero sections,
the \emph{self-normalized mass measure} are given by
\begin{equation}\label{e:mass-intro}
\cM_{s_k}:=\frac{|s_k|_{h^k}^2}{N_k}\,d\mu,
\qquad
\widehat{\cM}_{s_k}:=\frac{|s_k|_{h^k}^2}{\|s_k\|_k^2}\,d\mu.
\end{equation}
These are finite positive Borel measures on $K$, considered with the
weak topology. On $\{s_k=0\}$, we set $\widehat{\cM}_{s_k}=\mu$.
This convention affects only finitely many degrees almost surely,
since $\P(s_k=0)=\P(\xi=0)^{N_k}$ and $N_k\asymp k^n$.

Observe that orthonormality gives 
\begin{equation}\label{e:mass-normalization-relation}
r_k:=\cM_{s_k}(K)
=\frac1{N_k}\sum_{j=1}^{N_k}|\xi_{j,k}|^2,
\qquad
\cM_{s_k}=r_k\widehat{\cM}_{s_k}.
\end{equation}
In particular, the two normalizations are related by the random total mass $r_k$.

For a finite Borel measure $\rho$ on $K$ and
$f\in \mathcal C ^0 (K, \mathbb R)$, we use the notation
$\rho(f):=\int_K f\,d\rho$.
We refer to $\cM_{s_k}(f)$ and
$\widehat{\cM}_{s_k}(f)$ as the corresponding
mass statistics.

\subsection{Sharp mass equidistribution}

The first question is what moment assumptions on the coefficients
are needed for mass equidistribution. Bayraktar~\cite{Bayraktar2020}
established almost-sure mass equidistribution for random polynomials
with sub-Gaussian coefficients, using Hanson--Wright concentration,
and obtained a corresponding result for ample line bundles with
smooth metrics and an ambient volume form. Here we work with big line bundles equipped with
continuous metrics and Bernstein--Markov scalar products, replacing
the sub-Gaussian assumption by substantially weaker moment conditions.
We obtain necessary and sufficient conditions for convergence in
probability and, under pairwise independence of coefficient rows,
for almost-sure convergence. These results hold for every prescribed
sequence of orthonormal bases.

\begin{Theorem}
\label{T:fixed-basis-dimension}
Under the preceding assumptions, the following hold.
\begin{enumerate}[label=\textup{(\roman*)}]
\item For any coupling between rows,
\[
\cM_{s_k}\xrightarrow{\P}\mu_{\rm eq}
\quad\Longleftrightarrow\quad
\E|\xi|^2=1.
\]
\item If the coefficient vectors
$\xi_k=(\xi_{1,k},\ldots,\xi_{N_k,k})^T$ are pairwise independent,
\[
\cM_{s_k}\xrightarrow{\rm a.s.}\mu_{\rm eq}
\quad\Longleftrightarrow\quad
\E|\xi|^2=1\ \text{and}\ \E|\xi|^{2+2/n}<\infty.
\]
\end{enumerate}
The moment conditions in \textup{(ii)} are sufficient under arbitrary
coupling between rows, and the sufficient conditions in
both parts also imply the corresponding
convergence of $\widehat{\cM}_{s_k}$.
\end{Theorem}

The exponent $2+2/n$ arises from the dimension growth
$N_k\asymp k^n$. Indeed, the total mass is an average
of $N_k$ independent copies of $|\xi|^2$, and the
probabilities of its deviations from the mean are
summable precisely when $\E|\xi|^{2+2/n}<\infty$.
This is closely related to the classical complete-convergence
results of Hsu--Robbins~\cite{HsuRobbins1947} and
Erd\H{o}s~\cite{Erdos1949}, and to results along
subsequences due to Gut~\cite{Gut1985}.
Our theorem shows that the same moment condition also
suffices for almost-sure convergence of the mass measures,
for every prescribed sequence of orthonormal bases.
The proof reduces the mass statistics to Hermitian
quadratic forms with uniformly bounded operator norms,
for which we establish the required summable deviation
estimate in Appendix~\ref{S:quadratic-moment}.

Pairwise independence between coefficient rows is used
only to establish necessity of the higher moment condition
in part~\textup{(ii)}. Indeed, for nested rows drawn
from a single infinite iid sequence, the higher moment
condition can be dropped entirely, and $\E|\xi|^2=1$
already suffices for almost-sure mass equidistribution;
see Appendix~\ref{S:row-coupling} for details.

The proof of the theorem
is given in Section~\ref{S:Mass} and the corresponding results for self-normalized mass are
established in Corollaries~\ref{C:self-normalized-prob}
and~\ref{C:self-normalized-as}, where the unit-variance
assumption is relaxed to a finite positive second moment.

\subsection{From mass to zeros}

Mass and zeros reflect different aspects of a holomorphic section.
The former is determined by its
squared pointwise norm, whereas the latter is governed by the
logarithm of that norm through the Poincar\'e--Lelong formula.
It is therefore natural to ask when mass equidistribution
forces convergence of the zero currents. Such an implication is
not automatic, since convergence of the mass measures does
not directly control the logarithmic potentials. We show that mass equidistribution implies convergence of zero currents for big line bundles under an additional condition on the equilibrium measure.

In this part, we assume further that $X$ is K\"ahler and that
$\mu_{\rm eq}(\operatorname{Int}_X(\supp\mu_{\rm eq}))=1$.
Choose a smooth reference metric $h_0$ on
$L$, write $h=h_0e^{-\phi}$, and set $\theta=c_1(L,h_0)$.
The equilibrium potential $P_K\phi$ is defined in
Section~\ref{S:Geometry}.

\begin{Theorem}
Under these assumptions, let
$s_k\in H^0(X,L^k)\setminus\{0\}$ be a deterministic sequence.
If $\cM_{s_k}\to\mu_{\rm eq}$ weakly, then
\[
\frac1k\log|s_k|_{h_0^k}^2
\longrightarrow P_K\phi
\qquad\text{in }L^1(X),
\]
and
\[
\frac1k[Z_{s_k}]
\longrightarrow\theta+dd^cP_K\phi
\]
weakly as currents. If instead
$\widehat{\cM}_{s_k}\to\mu_{\rm eq}$, the same
conclusion for zero currents holds.
\end{Theorem}

Here $[Z_{s_k}]$ denotes the current of integration over the zero
divisor of $s_k$. The key point is that mass convergence identifies the
limiting logarithmic potential on the open set
$\operatorname{Int}_X(\supp\mu_{\rm eq})$. Locality of the
non-pluripolar Monge--Amp\`ere measure then gives full mass for
the limiting potential, and the domination principle extends
the identification to all of $X$. The proof is given in
Theorem~\ref{T:deterministic-mass-zero}. An example with equilibrium
measure supported on a circle shows that the assumption on
$\mu_{\rm eq}$ cannot in general be omitted.

This extends the mass-to-zero approach of Shiffman--Zelditch
\cite{ShiffmanZelditch1999} and Zelditch~\cite{Zelditch2018},
developed for weighted polynomials by
Bayraktar~\cite[Theorem~2.4]{Bayraktar2020}.
Zero equidistribution for big line bundles with continuous metrics
was also established by Bayraktar--Coman--Marinescu--Nguy\^en
\cite{BCMN2024}, using Fubini--Study currents and inner products
induced by an ambient volume form, under logarithmic-moment conditions
on coefficient functionals. Their results apply to more general
complex spaces. Our approach instead
gives a deterministic mass-to-zero implication for
Bernstein--Markov scalar products. Combined with the mass
equidistribution theorem, it yields convergence of zero currents
in probability and almost surely under the corresponding moment
assumptions, including for discrete coefficient laws such as
symmetric Bernoulli variables, which do not satisfy the
logarithmic-moment condition of~\cite{BCMN2024}.
We discuss this distinction further in Section~\ref{S:MassZeros}.

\subsection{Mass fluctuations and orthonormal bases}

Mass equidistribution describes the limiting behavior of random
sections, but does not capture the fluctuations of their mass
statistics. It is therefore natural to ask whether these
fluctuations satisfy a central limit theorem and how their
behavior depends on the coefficient distribution and the choice
of orthonormal basis.

In~\cite{BayraktarBojnik2026}, we establish a central limit
theorem for Gaussian mass statistics associated with sequences
of positive line bundles equipped with $\mathscr C^3$ metrics,
using the Sodin--Tsirelson theorem for nonlinear functionals
of Gaussian processes together with Bergman-kernel estimates.
Here, we work with big line bundles equipped with continuous
metrics and utilize a different approach. Namely, our approach treats 
general centered coefficient laws with finite fourth moments, 
without symmetry assumptions, and requires only weak 
Bergman-measure asymptotics.

Unlike the standard complex Gaussian ensemble, random sections
with general iid coefficients need not be invariant under
unitary changes of basis. In particular, the variance of a
mass statistic may depend on the basis; see
Lemma~\ref{P:mass-covariance}. Since the basis-dependent
contributions need not converge as $k\to\infty$, we establish
a central limit theorem using the exact finite-degree
variance, under a suitable nondegeneracy condition.

\begin{Theorem}[CLT]
\label{T:intro-exact-variance-CLT}
Assume $\E|\xi|^2=1$ and
$\nu_4:=\E|\xi|^4<\infty$.
Fix $f\in \mathcal C^0(K,\R)$. If
\begin{equation}\label{e:CLT-variance-lower-intro}
\liminf_{k\to\infty}
N_k\Var\bigl(\cM_{s_k}(f)\bigr)>0,
\end{equation}
then
\[
\frac{\cM_{s_k}(f)-\E\cM_{s_k}(f)}
{\sqrt{\Var(\cM_{s_k}(f))}}
\xrightarrow{\rm d}\mathcal N(0,1).
\]
In particular, condition~\eqref{e:CLT-variance-lower-intro}
holds whenever
\begin{equation}\label{e:CLT-nondegeneracy-intro}
|\E\xi^2|<1,\qquad
\nu_4>1,\qquad
\mu_{\rm eq}(f^2)>0.
\end{equation}
\end{Theorem}

The proof begins by expressing the mass statistics as
quadratic forms in the random coefficients and then
combines the quadratic-form central limit theorem of
Bai--Yao~\cite{BaiYao2008,BaiYao2021} with Toeplitz trace
asymptotics obtained from the Bergman-measure results of
Berman--Boucksom--Witt Nystr\"om~\cite{BBWN2011} and
Finski~\cite{Finski2026}. The covariance formula shows that
when $\E\xi^2=0$ and $\nu_4=2$, the variance agrees with
that of the standard complex Gaussian ensemble, even for
non-Gaussian coefficients. In general, however, the choice
of basis can have a substantial effect.
Example~\ref{Ex:projective-bases} shows that, for coefficients
uniformly distributed on the unit circle, a mass statistic
has identically vanishing fluctuations in the monomial basis
on $\mathbf P^1$, but nondegenerate Gaussian fluctuations
after a suitable unitary change of basis.

We also establish a central limit theorem for self-normalized
mass statistics in Theorem~\ref{T:self-normalized-CLT}.
Although these central limit theorems apply to arbitrary
bases under suitable nondegeneracy conditions, their
normalized variances need not converge. This motivates
asking whether typical choices of bases yield explicit
limiting variances. In the Haar random-basis framework
of Zelditch~\cite{ZelditchRandomBases2014}, we show that
the basis-dependent covariance terms have deterministic
limits for almost every sequence of bases.
The resulting limiting variances at the $\sqrt{N_k}$
scale are
\[
\begin{aligned}
\sigma^2(f)
&=\mu_{\rm eq}(f^2)
  +(\nu_4-2)\mu_{\rm eq}(f)^2,\\
\widehat{\sigma}^2(f)
&=\mu_{\rm eq}(f^2)
  -\mu_{\rm eq}(f)^2,
\end{aligned}
\]
for dimension-normalized and self-normalized mass statistics,
respectively. In particular, the self-normalized limiting
variance depends only on the equilibrium measure and the
test function, independently of the coefficient law, giving rise to a universality phenomenon for self-normalized fluctuations. The precise statements are given in
Proposition~\ref{T:Haar-basis-CLT}. 

These results are closely related to quantum ergodicity of
random orthonormal bases; see
Zelditch~\cite{ZelditchRandomBases2014,Zelditch2018}.
Indeed, the first profile limit in
Proposition~\ref{T:Haar-basis-CLT} implies that, for
Haar-almost every sequence of orthonormal bases
$\{S_{j,k}\}_{j=1}^{N_k}$,
\[
\frac{1}{N_k}\sum_{j=1}^{N_k}
\left|
\int_K f\,|S_{j,k}|_{h^k}^2\,d\mu
-\mu_{\rm eq}(f)
\right|^2
\longrightarrow 0
\]
for every $f\in\mathcal C^0(K,\mathbb R)$.
Thus, these bases are quantum ergodic with equilibrium state
$f\mapsto\mu_{\rm eq}(f)$. Our central limit theorems complement
this equidistribution property by describing the fluctuations
of random sections with iid coefficients after almost every
such sequence of bases has been fixed.

\medskip
\noindent\textbf{Organization.}
In Section~\ref{S:Geometry}  we recall the geometric background
and Toeplitz trace asymptotics. The mass equidistribution
and mass-to-zero results are proved in
Sections~\ref{S:Mass} and~\ref{S:MassZeros}, respectively.
Section~\ref{S:CLT} contains the fluctuation theorems,
and the auxiliary probabilistic results appear in the appendices.

\section{Equilibrium measures and Toeplitz operators}
\label{S:Geometry}

In this section, we recall the equilibrium potential and measure
associated with $(K,\phi)$, together with the Bergman kernel
and Toeplitz trace asymptotics needed in our proofs.
The diagonal Bergman limit identifies the limiting mass measure,
while the Toeplitz trace asymptotics enter the fluctuation theorems.

\subsection{Equilibrium potentials and measures}

We retain the geometric setting of the introduction.
Throughout, $dd^c=\frac{\imat}{2\pi}\partial\bar\partial$.
Convergence in $L^1(X)$ is understood with respect to a fixed
smooth positive volume form on $X$, whose choice is immaterial
since $X$ is compact.

Recall that the volume of $L$ is defined by
\begin{equation}\label{e:volume}
\vol(L):=\limsup_{k\to\infty}
\frac{n!}{k^n}h^0(X,L^k)>0.
\end{equation}
For a big line bundle, the limsup is in fact a limit, giving
\eqref{e:dimension-asymptotic-intro}.
In particular, there exist constants $c,C>0$ such that
\begin{equation}\label{e:N-comparable}
ck^n\leq N_k\leq Ck^n
\end{equation}
for all sufficiently large $k$; see
\cite[Section~2.2]{BermanBoucksom2010}.

Fix a smooth Hermitian metric $h_0$ on $L$ and set
$\theta=c_1(L,h_0)$. Write $h=h_0e^{-\phi}$, where
$\phi\in\mathcal{C}^0(X,\mathbb R)$, with the convention
$|s|_h^2=e^{-\phi}|s|_{h_0}^2$.
If $e$ is a local holomorphic frame with
$|e|_{h_0}^2=e^{-\rho}$, then
$\theta=dd^c\rho$ and $|e|_h^2=e^{-(\rho+\phi)}$.

We denote the space of $\theta$-plurisubharmonic functions by
\[
\operatorname{PSH}(X,\theta)
:=\bigl\{u\in L^1(X):\theta+dd^cu\geq0\bigr\},
\]
where the inequality is understood in the sense of currents
and each function is identified with its unique
$\theta$-psh representative. Equivalently, $\rho+u$
is plurisubharmonic in every local frame.
We say that $u\in\operatorname{PSH}(X,\theta)$ has
\emph{minimal singularities} if, for every
$v\in\operatorname{PSH}(X,\theta)$, there exists
$C_v\in\mathbb R$ such that $v\leq u+C_v$ on $X$.
When $L$ is merely big, potentials with minimal
singularities need not be locally bounded; see
\cite[Section~1.4]{BEGZ2010}.

Recall that a subset of $X$ is pluripolar if it is
locally contained in the $-\infty$ locus of a
plurisubharmonic function that is not identically
$-\infty$. The weighted equilibrium potential
associated with $(K,\phi)$ is defined as
\begin{equation}\label{e:equilibrium-weight}
P_K\phi:=
\left(
\sup\left\{
u\in\operatorname{PSH}(X,\theta):
u\leq\phi\ \text{on }K
\right\}
\right)^*,
\end{equation}
where $*$ denotes upper semicontinuous regularization.
Since $K$ is non-pluripolar, $P_K\phi$ belongs to
$\operatorname{PSH}(X,\theta)$ and has minimal
singularities; see
\cite[Section~2.4]{BermanBoucksom2010}.

The associated equilibrium current is
$T_{\rm eq}:=\theta+dd^cP_K\phi$. Since $L$ is
merely big, this current may have singular potentials,
and its top exterior power is understood in the
non-pluripolar sense. The equilibrium measure is
defined by
\begin{equation}\label{e:equilibrium-measure}
\mu_{\rm eq}:=
\frac{1}{\vol(L)}
\big\langle T_{\rm eq}^{\,n}\big\rangle
=
\frac{1}{\vol(L)}
\big\langle(\theta+dd^cP_K\phi)^n\big\rangle.
\end{equation}
The brackets denote the non-pluripolar
Monge--Amp\`ere product. It agrees with the
Bedford--Taylor product wherever the potential
is locally bounded and does not charge
pluripolar sets. For its construction and
basic properties, see
\cite[Section~1]{BEGZ2010}.

Since $P_K\phi$ has minimal singularities,
\cite[Theorem~2.2]{BermanBoucksom2010} yields
\[
\int_X\big\langle T_{\rm eq}^{\,n}\big\rangle
=\vol(L).
\]
Moreover, \cite[Proposition~2.10]{BermanBoucksom2010}
shows that $\mu_{\rm eq}$ is concentrated on $K$ and
that $P_K\phi=\phi$ holds $\mu_{\rm eq}$-almost
everywhere. Thus $\mu_{\rm eq}$ is a probability
measure supported on $K$. In particular, no weighted
regularity assumption on $(K,\phi)$ is required.

In the special case $K=X$, with $h$ smooth and
positively curved, one has $P_X\phi=\phi$, and hence
\[
T_{\rm eq}=c_1(L,h),
\qquad
\mu_{\rm eq}=\frac{c_1(L,h)^n}{\vol(L)}.
\]

\subsection{Bernstein--Markov measures and Bergman kernel asymptotics}
\label{S:Bergman-kernels}

We remark that the inner product~\eqref{e:L2-intro} is positive definite.
Indeed, if $\|s\|_k=0$, the continuity of
$|s|_{h^k}$ implies that $s$ vanishes on
$K=\supp\mu$. Since $K$ is non-pluripolar,
whereas the zero set of a nonzero holomorphic section
on the connected manifold $X$ is pluripolar, it follows
that $s\equiv0$. We write $\|s\|_k^2=(s,s)_k$.

Recall that $\mu$ is a Bernstein--Markov measure
for $(K,h)$, meaning that for every $\varepsilon>0$,
there exists $C_\varepsilon>0$ such that
\begin{equation}\label{e:BM}
\sup_K|s|_{h^k}
\leq C_\varepsilon e^{\varepsilon k}\|s\|_k,
\qquad
s\in H^0(X,L^k),\quad k\geq1.
\end{equation}
Since $\mu$ is a probability measure, this property
ensures that the supremum and $L^2(\mu)$ norms
differ by at most subexponential factors.
We refer to~\cite{BloomLevenbergSurvey,Finski2026}
for examples and further references.

We note that the non-pluripolarity assumption concerns the
support $K$, the measure $\mu$ itself is not required
to vanish on pluripolar sets.

For all sufficiently large $k$, let
$\{S_{1,k},\ldots,S_{N_k,k}\}$ be an orthonormal
basis of the Hilbert space
$(\cH_k,(\cdot,\cdot)_k)$.
The associated Bergman kernel is then given by
\begin{equation}\label{e:Bergman-kernel}
B_k(x,y):=
\sum_{j=1}^{N_k}
S_{j,k}(x)\otimes S_{j,k}(y)^*,
\end{equation}
where $S_{j,k}(y)^*\in(L_y^k)^*$ denotes the metric dual
of $S_{j,k}(y)$, defined by
$S_{j,k}(y)^*(v)=h_y^k(v,S_{j,k}(y))$
for $v\in L_y^k$.
The kernel is independent of the chosen
orthonormal basis, and its diagonal density is denoted by
\begin{equation}\label{e:Bergman-diagonal}
B_k(x):=
\sum_{j=1}^{N_k}|S_{j,k}(x)|_{h^k}^2.
\end{equation}
Orthonormality gives the identity
\begin{equation}\label{e:reproducing-norm}
\int_K
|B_k(x,y)|_{h_x^k\otimes(h_y^k)^*}^2
\,d\mu(y)=B_k(x).
\end{equation}

Following~\cite{Finski2026}, we introduce the
normalized diagonal and off-diagonal Bergman measures
\begin{align}
\mu_k^{\mathrm{Berg},\Delta}
&:=\frac{1}{N_k}B_k(x)\,d\mu(x),
\label{e:beta}\\
\mu_k^{\mathrm{Berg}}
&:=\frac{1}{N_k}
|B_k(x,y)|_{h_x^k\otimes(h_y^k)^*}^2
\,d\mu(x)d\mu(y).
\label{e:Gamma}
\end{align}
These are probability measures on $X$ and
$X\times X$, respectively. Indeed,
$\int_K B_k\,d\mu=N_k$, and
\eqref{e:reproducing-norm} shows that both
marginals of $\mu_k^{\mathrm{Berg}}$ coincide
with $\mu_k^{\mathrm{Berg},\Delta}$.

The asymptotic behavior of these measures is described
by the diagonal Bergman measure convergence theorem of
Berman--Boucksom--Witt Nystr\"om
\cite[Theorem~B]{BBWN2011} and the off-diagonal
convergence theorem of Finski
\cite[Theorem~1.3]{Finski2026}.
The latter combines the diagonal convergence with
the asymptotic concentration of the off-diagonal
Bergman measures along the diagonal.
We recall both results below.

\begin{Theorem}[Bergman measure limits]
\label{T:Bergman-limits}
Under the preceding assumptions,
\begin{align}
\mu_k^{\mathrm{Berg},\Delta}
&\longrightarrow\mu_{\rm eq},
\label{e:diagonal-limit}\\
\mu_k^{\mathrm{Berg}}
&\longrightarrow\Delta_*\mu_{\rm eq},
\label{e:offdiagonal-limit}
\end{align}
weakly as measures on $X$ and $X\times X$,
respectively, where $\Delta:X\to X\times X$,
$\Delta(x)=(x,x)$, is the diagonal embedding.
\end{Theorem}

Here $\Delta_*\mu_{\rm eq}$ denotes the pushforward
of $\mu_{\rm eq}$ under $\Delta$, characterized by
\[
\int_{X\times X}\Psi\,d(\Delta_*\mu_{\rm eq})
=
\int_X\Psi(x,x)\,d\mu_{\rm eq}(x),
\qquad \Psi\in \mathcal C^0(X\times X).
\]
In particular, for every $\Psi\in \mathcal C^0(X\times X)$,
\begin{equation}\label{e:Gamma-test}
\begin{aligned}
&\frac{1}{N_k}\int_{K\times K}
\Psi(x,y)
|B_k(x,y)|_{h_x^k\otimes(h_y^k)^*}^2
\,d\mu(x)d\mu(y)\longrightarrow
\int_K\Psi(x,x)\,d\mu_{\rm eq}(x).
\end{aligned}
\end{equation}
Since $K\times K$ is closed in $X\times X$,
the Tietze extension theorem shows that
\eqref{e:Gamma-test} remains valid for every
$\Psi\in \mathcal C^0(K\times K)$. We shall use this
convergence with $\Psi(x,y)=f(x)g(y)$, together
with \eqref{e:diagonal-limit}, in the
Toeplitz trace asymptotics below.

\subsection{Toeplitz operators and trace asymptotics}

Let $\Pi_k:L^2(K,L^k;h^k,\mu)\to\cH_k$ denote the
orthogonal projection. For $f\in L^\infty(K)$, the Toeplitz
operator with symbol $f$ is defined by
\begin{equation}\label{e:Toeplitz-def}
T_k(f):=\Pi_k\circ M_f\big|_{\cH_k},
\end{equation}
where $M_f$ is multiplication by $f$. Equivalently, $T_k(f)$
is characterized by
\begin{equation}\label{e:Toeplitz-form}
(T_k(f)s,t)_k
=\int_K f\,h^k(s,t)\,d\mu,
\qquad s,t\in\cH_k.
\end{equation}
Indeed, the right-hand side is a bounded sesquilinear form,
so the Riesz representation theorem applies. Moreover,
$T_k(f)^*=T_k(\overline f)$; in particular, $T_k(f)$ is
self-adjoint when $f$ is real-valued.

Fix the orthonormal basis
$\{S_{1,k},\ldots,S_{N_k,k}\}$
used in~\eqref{e:random-section-intro}. The matrix of $T_k(f)$
in this basis is $A_k(f)=\bigl(a_{ij}^{(k)}(f)\bigr)_{i,j=1}^{N_k}$, where
\begin{equation}\label{e:A-entry}
\begin{aligned}
a_{ij}^{(k)}(f)
&=(T_k(f)S_{j,k},S_{i,k})_k
=\int_K f(x)h^k(S_{j,k}(x),S_{i,k}(x))\,d\mu(x).
\end{aligned}
\end{equation}
For real-valued $f$, the matrix $A_k(f)$ is Hermitian.
The trace is independent of the orthonormal basis, and
$\Tr T_k(f)=\Tr A_k(f)=\sum_{j=1}^{N_k}a_{jj}^{(k)}(f)$.

We write $\|A\|_{\rm op}$ for the operator norm and
$\|A\|_{\rm HS}^2=\Tr(A^*A)=\sum_{i,j}|a_{ij}|^2$
for the unnormalized Hilbert--Schmidt square.
Orthonormality identifies these matrix norms with the
corresponding norms of the Toeplitz operator.
For real-valued $f$, Hermitian symmetry gives
\[
\Tr\bigl(A_k(f)^2\bigr)
=\norm{A_k(f)}_{\rm HS}^2
=\sum_{j=1}^{N_k}a_{jj}^{(k)}(f)^2
+\sum_{i\ne j}|a_{ij}^{(k)}(f)|^2.
\]
The trace on the left is invariant under unitary changes of
basis, although the two sums on the right need not be
separately invariant.

Since $\Pi_k$ is an orthogonal projection and
$\norm{M_f}_{\rm op}\leq\norm{f}_\infty$, we obtain
\begin{equation}\label{e:Toeplitz-op-bound}
\norm{A_k(f)}_{\rm op}
=\norm{T_k(f)}_{\rm op}\leq\norm{f}_\infty.
\end{equation}
The inequality $\norm{A}_{\rm HS}^2\leq
N\norm{A}_{\rm op}^2$ for $A\in\C^{N\times N}$
therefore yields, for real-valued $f$,
\[
\Tr\bigl(A_k(f)^2\bigr)
=\norm{A_k(f)}_{\rm HS}^2
\leq N_k\norm{A_k(f)}_{\rm op}^2
\leq N_k\norm{f}_\infty^2.
\]
Also, for every $f\in L^\infty(K)$ and every $j$,
\[
|a_{jj}^{(k)}(f)|
=|(T_k(f)S_{j,k},S_{j,k})_k|
\leq\norm{T_k(f)}_{\rm op}\leq\norm{f}_\infty.
\]

\begin{Proposition}[Toeplitz trace asymptotics]
\label{P:trace-identities}
For every $f,g\in \mathcal C^0(K,\R)$,
\begin{align}
\frac{1}{N_k}\Tr A_k(f)
&\longrightarrow\int_K f\,d\mu_{\rm eq},
\label{e:first-trace-limit}\\
\frac{1}{N_k}\Tr\bigl(A_k(f)A_k(g)\bigr)
&\longrightarrow\int_K fg\,d\mu_{\rm eq}.
\label{e:second-trace-limit}
\end{align}
In particular,
\begin{equation}\label{e:HS-limit}
\frac{1}{N_k}\Tr\bigl(A_k(f)^2\bigr)
\longrightarrow\int_K f^2\,d\mu_{\rm eq}.
\end{equation}
\end{Proposition}

\begin{proof}
By \eqref{e:A-entry} and \eqref{e:Bergman-diagonal},
\[
\begin{aligned}
\frac{1}{N_k}\Tr A_k(f)
&=\frac{1}{N_k}\sum_{j=1}^{N_k}a_{jj}^{(k)}(f)=\frac{1}{N_k}\int_K f(x)B_k(x)\,d\mu(x)
=\int_K f\,d\mu_k^{\mathrm{Berg},\Delta}.
\end{aligned}
\]
The first limit follows from
$\mu_k^{\mathrm{Berg},\Delta}\to\mu_{\rm eq}$ in
Theorem~\ref{T:Bergman-limits}.

For the second limit, Fubini's theorem and
\eqref{e:A-entry} give
\[
\begin{aligned}
\frac{1}{N_k}\Tr\bigl(A_k(f)A_k(g)\bigr)
&=\frac{1}{N_k}\sum_{i,j=1}^{N_k}
a_{ij}^{(k)}(f)a_{ji}^{(k)}(g)\\
&=\frac{1}{N_k}\int_{K\times K}f(x)g(y)
|B_k(x,y)|_{h_x^k\otimes(h_y^k)^*}^2
\,d\mu(x)d\mu(y)\\
&=\int_{K\times K}f(x)g(y)\,
d\mu_k^{\mathrm{Berg}}(x,y).
\end{aligned}
\]
The convergence $\mu_k^{\mathrm{Berg}}\to
\Delta_*\mu_{\rm eq}$ in Theorem~\ref{T:Bergman-limits}
yields \eqref{e:second-trace-limit}. Taking $g=f$
then gives \eqref{e:HS-limit}.
\end{proof}

\section{Mass equidistribution}\label{S:Mass}

In this section, we prove Theorem~\ref{T:fixed-basis-dimension}.
The coefficients are independent within each row, while
independence between rows is used only in the proof of the
almost-sure converse.

Before turning to the proofs, we recall a convenient
metrization of the weak topology of finite positive measures, which will be useful to pass from mass statistics against individual test functions to weak convergence of the random measures.

Choose a uniformly dense sequence
$(f_\ell)_{\ell\geq1}$ in
$\mathcal C^0(K,\R)$ with $f_1\equiv1$.
For finite positive measures $\rho,\eta$ on $K$,
set
\begin{equation}\label{e:weak-metric}
d_{\rm w}(\rho,\eta)
=\sum_{\ell\geq1}2^{-\ell}
\frac{|\rho(f_\ell)-\eta(f_\ell)|}
{1+|\rho(f_\ell)-\eta(f_\ell)|}.
\end{equation}
This metric induces the weak topology. Indeed,
convergence against the dense family implies
convergence of the total masses, since $f_1=1$,
and uniform approximation then gives convergence
against every continuous function. Conversely,
weak convergence implies convergence in this metric
by dominated convergence for the series.

\subsection{Convergence in probability}
\label{S:dimension-probability}

We begin with the quadratic-form representation of the mass
statistics. Since Hermitian products are linear in the first
variable, we have
\[
|s_k(x)|_{h^k}^2
=\sum_{i,j=1}^{N_k}\ov\xi_{i,k}\xi_{j,k}\,
h^k(S_{j,k}(x),S_{i,k}(x)).
\]
Integrating against $f\,d\mu$, for $f\in\mathcal C^0(K,\R)$, and using
\eqref{e:A-entry}, we obtain
\begin{equation}\label{e:mass-stat}
\cM_{s_k}(f)=\frac1{N_k}\xi_k^*A_k(f)\xi_k.
\end{equation}
If $\E|\xi|^2=1$, independence and the standing centering
assumption give $\E(\ov\xi_{i,k}\xi_{j,k})=\delta_{ij}$.
Consequently, we have
\begin{equation}\label{e:expectation}
\E\cM_{s_k}(f)
=\frac1{N_k}\Tr A_k(f)
=\mu_k^{\mathrm{Berg},\Delta}(f).
\end{equation}
By \eqref{e:diagonal-limit}, these expectations converge
to $\mu_{\rm eq}(f)$. To prove convergence in probability,
it remains to control the fluctuations of the mass statistics
around their expectations.

\begin{proof}[Proof of Theorem~\ref{T:fixed-basis-dimension}\textup{(i)}]
\emph{Sufficiency.}
Assume $\E|\xi|^2=1$. We first establish convergence
in mean of the mass statistics. Fix
$f\in\mathcal C^0(K,\R)$ and write
$A_k=A_k(f)=(a_{ij}^{(k)})$. By
\eqref{e:mass-stat} and \eqref{e:expectation},
\[
\cM_{s_k}(f)-\E\cM_{s_k}(f)=D_k+O_k,
\]
where
\[
D_k=\frac1{N_k}\sum_j a_{jj}^{(k)}
\bigl(|\xi_{j,k}|^2-1\bigr),
\qquad
O_k=\frac1{N_k}\sum_{i\ne j}
a_{ij}^{(k)}\ov\xi_{i,k}\xi_{j,k}.
\]
The two terms require different estimates. The diagonal term
contains $|\xi|^2-1$, whose first moment is finite but whose variance
may be infinite; we therefore truncate it. In the off-diagonal term,
independence and centering reduce the second-moment calculation to
pairings involving at most two factors of each coefficient.

For the diagonal term, set $Y=|\xi|^2-1$ and
$Y_{j,k}=|\xi_{j,k}|^2-1$. The variables $Y_{j,k}$
have the common law of $Y$, which is centered
and integrable. Setting $b_{j,k}=a_{jj}^{(k)}/N_k$,
the operator-norm bound \eqref{e:Toeplitz-op-bound}
gives $\max_j|b_{j,k}|\leq\|f\|_\infty/N_k\to0$
and $\sum_j|b_{j,k}|\leq\|f\|_\infty$.
With this notation, $D_k=\sum_jb_{j,k}Y_{j,k}$.

Fix $R>0$ and truncate and recenter $Y$ by setting
\[
Y^{(R)}
=Y\1_{\{|Y|\leq R\}}
-\E\bigl(Y\1_{\{|Y|\leq R\}}\bigr),
\qquad
Z^{(R)}=Y-Y^{(R)}.
\]
Define $Y_{j,k}^{(R)}$ and $Z_{j,k}^{(R)}$
analogously from $Y_{j,k}$. Then
$D_k=\sum_jb_{j,k}Y_{j,k}^{(R)}
+\sum_jb_{j,k}Z_{j,k}^{(R)}$.
The variables $Y_{j,k}^{(R)}$ are centered, independent
within each row, and bounded in absolute value by $2R$.
Independence therefore yields
\begin{align*}
\E\left|\sum_jb_{j,k}Y_{j,k}^{(R)}\right|^2
&=\sum_jb_{j,k}^2\E|Y^{(R)}|^2
\leq4R^2\left(\max_j|b_{j,k}|\right)
\sum_j|b_{j,k}|
\longrightarrow0.
\end{align*}
For each fixed $R$, Cauchy--Schwarz then shows that
the first sum converges to zero in $L^1$.

Since $\E Y=0$, we have
$Z^{(R)}=Y\1_{\{|Y|>R\}}-
\E(Y\1_{\{|Y|>R\}})$, and hence
$\E|Z^{(R)}|
\leq2\E(|Y|\1_{\{|Y|>R\}})$.
The remainder satisfies
\[
\begin{aligned}
\E\left|\sum_jb_{j,k}Z_{j,k}^{(R)}\right|
&\leq\sum_j|b_{j,k}|\E|Z^{(R)}|
\leq2\|f\|_\infty
\E\bigl(|Y|\1_{\{|Y|>R\}}\bigr).
\end{aligned}
\]
Combining the two estimates and taking the upper limit
as $k\to\infty$ gives
\[
\limsup_{k\to\infty}\E|D_k|
\leq2\|f\|_\infty
\E\bigl(|Y|\1_{\{|Y|>R\}}\bigr).
\]
The right-hand side tends to zero as $R\to\infty$
by integrability of $Y$. Thus $D_k\to0$ in $L^1$.

We now turn to the off-diagonal term. Put $\beta=\E\xi^2$,
so that $|\beta|\leq\E|\xi|^2=1$. Expanding its second
moment gives
\[
N_k^2\E|O_k|^2
=\sum_{i\ne j}\sum_{p\ne q}
a_{ij}^{(k)}\ov{a_{pq}^{(k)}}
\E\bigl(\ov\xi_{i,k}\xi_{j,k}
\xi_{p,k}\ov\xi_{q,k}\bigr).
\]
By independence and centering, the expectation vanishes
whenever an index occurs only once. Since $i\ne j$ and
$p\ne q$, the only potentially nonvanishing terms
correspond to $(p,q)=(i,j)$ and $(p,q)=(j,i)$, with
expectations $1$ and $|\beta|^2$, respectively.
Using the Hermitian symmetry
$a_{ji}^{(k)}=\ov{a_{ij}^{(k)}}$, we obtain
\begin{equation}\label{e:offdiag-L2}
\begin{aligned}
\E|O_k|^2
&=\frac1{N_k^2}
\left(
\sum_{i\ne j}|a_{ij}^{(k)}|^2
+|\beta|^2\sum_{i\ne j}(a_{ij}^{(k)})^2
\right)\leq\frac{2\Tr(A_k^2)}{N_k^2}
\leq\frac{2\|f\|_\infty^2}{N_k}.
\end{aligned}
\end{equation}
The first inequality follows from
$\left|\sum_{i\ne j}(a_{ij}^{(k)})^2\right|
\leq\sum_{i\ne j}|a_{ij}^{(k)}|^2$
and the identity
$\Tr(A_k^2)=\sum_{i,j}|a_{ij}^{(k)}|^2$;
the second uses the operator-norm estimate
$\|A_k\|_{\rm op}\leq\|f\|_\infty$.

This calculation requires only second moments, since
each nonvanishing expectation involves at most two
factors of any given coefficient. As $N_k\to\infty$,
we have $O_k\to0$ in $L^2$, and hence in $L^1$.

Combining these estimates with \eqref{e:expectation}
and the diagonal Bergman-measure limit yields
\[
\begin{aligned}
\E\bigl|\cM_{s_k}(f)-\mu_{\rm eq}(f)\bigr|
&\leq\E|D_k|+\E|O_k|+
\bigl|\mu_k^{\mathrm{Berg},\Delta}(f)
-\mu_{\rm eq}(f)\bigr|
\longrightarrow0.
\end{aligned}
\]
Thus $\cM_{s_k}(f)\to\mu_{\rm eq}(f)$ in probability
for every $f\in \mathcal C^0(K,\R)$. Applying the same estimate
to the real and imaginary parts also gives $L^1$
convergence for every $f\in \mathcal C ^0 (K, \mathbb C)$.

To obtain weak convergence of the random measures,
we use the metric~\eqref{e:weak-metric}.
Fix $\varepsilon>0$ and choose $L$ such that
$\sum_{\ell>L}2^{-\ell}<\varepsilon/2$.
Since each $\cM_{s_k}(f_\ell)$ converges to
$\mu_{\rm eq}(f_\ell)$ in probability, a finite
union bound gives
\[
\begin{aligned}
&\P\bigl(d_{\rm w}(\cM_{s_k},\mu_{\rm eq})
>\varepsilon\bigr)\leq
\sum_{\ell=1}^L
\P\left(
|\cM_{s_k}(f_\ell)-\mu_{\rm eq}(f_\ell)|
>\frac{\varepsilon}{2}
\right)
\longrightarrow0.
\end{aligned}
\]
This proves weak convergence in probability. Note that no
independence between different rows has been used.

\smallskip\noindent

\emph{Necessity.}
Suppose that $\cM_{s_k}\to\mu_{\rm eq}$ weakly in
probability. Since $\mu_{\rm eq}$ is a probability measure,
testing against the constant function $1$ gives
\[
r_k=\cM_{s_k}(1)
=\frac1{N_k}\sum_{j=1}^{N_k}V_{j,k}
\longrightarrow1
\qquad\text{in probability},
\]
where $V=|\xi|^2$ and $V_{j,k}=|\xi_{j,k}|^2$.

We first show that $\E V<\infty$. Suppose, to the
contrary, that $\E V=\infty$. For $R>0$, set
$V^{(R)}=\min\{V,R\}$ and
$V_{j,k}^{(R)}=\min\{V_{j,k},R\}$. For each fixed
$R$, the variables $V_{j,k}^{(R)}$ are independent
and identically distributed within each row, with
$0\leq V_{j,k}^{(R)}\leq R$. Hence
\[
\begin{aligned}
\E\left|
\frac1{N_k}\sum_{j=1}^{N_k}V_{j,k}^{(R)}
-\E V^{(R)}
\right|^2
&=\frac{\Var(V^{(R)})}{N_k}\leq\frac{R^2}{N_k}
\longrightarrow0.
\end{aligned}
\]
Thus the truncated averages converge to $\E V^{(R)}$
in $L^2$, and therefore in probability.

Since $V$ is nonnegative, the monotone convergence
theorem gives $\E V^{(R)}\to\infty$ as $R\to\infty$.
Given $M>0$, choose $R$ such that $\E V^{(R)}>2M$.
As $V_{j,k}\geq V_{j,k}^{(R)}$, we obtain
\[
\begin{aligned}
\P(r_k>M)
&\geq
\P\left(
\frac1{N_k}\sum_{j=1}^{N_k}V_{j,k}^{(R)}>M
\right)
\longrightarrow1.
\end{aligned}
\]
Since $M$ is arbitrary, $r_k\to\infty$ in probability,
contradicting the assumed convergence to $1$.

Consequently, $\E V<\infty$. The ordinary weak law
of large numbers, applied to the independent and
identically distributed variables within each row,
now gives $r_k\to\E V$ in probability, since
$N_k\to\infty$. No independence between different
rows is required. By uniqueness of the limit in
probability, $\E V=1$, or equivalently
$\E|\xi|^2=1$.
\end{proof}

The sufficiency argument also shows that, if
$\E|\xi|^2=\sigma^2\in(0,\infty)$, then
$\cM_{s_k}\to\sigma^2\mu_{\rm eq}$ weakly in
probability. Indeed, applying the result to
$\xi/\sigma$ gives the conclusion, since the
corresponding mass measures are
$\sigma^{-2}\cM_{s_k}$.

We now turn to self-normalized mass measures.
Since self-normalization is invariant under rescaling,
the coefficients need not have unit second moment.

\begin{Corollary}
\label{C:self-normalized-prob}
Suppose that $0<\sigma^2:=\E|\xi|^2<\infty$.
Then $\widehat{\cM}_{s_k}\to\mu_{\rm eq}$ weakly
in probability, without any independence assumption
between coefficient rows.
\end{Corollary}

\begin{proof}
Set $t_k=s_k/\sigma$. Its coefficients are centered
and have unit second moment. By
Theorem~\ref{T:fixed-basis-dimension}\textup{(i)},
$\cM_{t_k}\to\mu_{\rm eq}$ weakly in probability.
In particular, $\cM_{t_k}(K)\to1$ in probability.

Whenever $s_k\ne0$, rescaling does not change the
self-normalized measure, and for every
$f\in\mathcal C^0(K,\R)$,
\[
\widehat{\cM}_{s_k}(f)
=\widehat{\cM}_{t_k}(f)
=\frac{\cM_{t_k}(f)}{\cM_{t_k}(K)}.
\]
Since the numerator converges in probability to
$\mu_{\rm eq}(f)$ and the denominator to $1$,
the ratio converges in probability to $\mu_{\rm eq}(f)$.
Using the metric~\eqref{e:weak-metric}, we obtain
weak convergence of the measures in probability.

The convention on zero sections causes no difficulty.
Indeed, writing $q=\P(\xi=0)<1$, independence within
each row gives $\P(s_k=0)=q^{N_k}\to0$.
\end{proof}

Unlike the dimension-normalized case, the self-normalized
mass measures have total mass one. Thus the argument
based on the random total mass does not yield necessary
moment conditions for convergence in probability or
almost surely. We therefore restrict ourselves to
sufficient conditions for self-normalized convergence.

\subsection{Almost-sure convergence}
\label{S:almost-sure-proof}

In this part, we prove almost-sure convergence result which is based on summability of the
deviation probabilities for the quadratic forms
representing the mass statistics. The diagonal terms
determine the moment threshold, while the off-diagonal
terms must be controlled without imposing additional
moment assumptions. The following lemma provides the
required estimate and its proof is given in
Appendix~\ref{S:quadratic-moment}.

\begin{Lemma}[Complete convergence of quadratic forms]
\label{P:complete-convergence}
Suppose $N_k\asymp k^n$, $n\geq1$, and the rows
$z_k=(\xi_{1,k},\ldots,\xi_{N_k,k})^T\in\C^{N_k}$
have independent coordinates with a common law satisfying
\[
\E\xi=0,\qquad
\E|\xi|^2=1,\qquad
\E|\xi|^{2+2/n}<\infty.
\]
For any sequence of deterministic Hermitian matrices
$A_k$ with $\sup_k\|A_k\|_{\rm op}<\infty$, and every
$\varepsilon>0$,
\begin{equation}\label{e:complete-convergence}
\sum_k\P\bigl(
|z_k^*A_kz_k-\Tr A_k|>\varepsilon N_k
\bigr)<\infty.
\end{equation}
No independence between rows is required.
\end{Lemma}

\begin{proof}[Proof of Theorem~\ref{T:fixed-basis-dimension}
\textup{(ii)}: Sufficiency]
Assume $\E|\xi|^2=1$ and $\E|\xi|^{2+2/n}<\infty$.
Fix $f\in\mathcal C^0(K,\R)$.
The matrices $A_k(f)$ are Hermitian, and
\eqref{e:Toeplitz-op-bound} gives
$\|A_k(f)\|_{\rm op}\leq\|f\|_\infty$.

Applying Lemma~\ref{P:complete-convergence} to $A_k(f)$
and using \eqref{e:mass-stat} and \eqref{e:expectation},
we obtain, for every $\varepsilon>0$,
\[
\sum_k\P\left(
\left|\cM_{s_k}(f)
-\mu_k^{\mathrm{Berg},\Delta}(f)\right|
>\varepsilon
\right)<\infty.
\]
By the first Borel--Cantelli lemma, for each
$\varepsilon>0$, the event
$\{|\cM_{s_k}(f)-\mu_k^{\mathrm{Berg},\Delta}(f)|
>\varepsilon\}$ occurs only finitely often almost
surely. Applying this with $\varepsilon=1/m$, $m\geq1$,
and intersecting the resulting probability-one events,
we conclude that
\[
\cM_{s_k}(f)
-\mu_k^{\mathrm{Berg},\Delta}(f)
\longrightarrow0
\qquad\text{almost surely}.
\]
Since $\mu_k^{\mathrm{Berg},\Delta}(f)
\to\mu_{\rm eq}(f)$ by \eqref{e:diagonal-limit},
it follows that $\cM_{s_k}(f)\to\mu_{\rm eq}(f)$
almost surely for every fixed $f\in\mathcal C^0(K,\R)$.

To pass to weak convergence of the random measures,
let $(f_\ell)_{\ell\geq1}$ be the dense family
defining the metric~\eqref{e:weak-metric}.
For each $\ell$, the preceding argument gives an
event $\Omega_\ell$ of probability one on which
$\cM_{s_k}(f_\ell)\to\mu_{\rm eq}(f_\ell)$.
Their intersection
$\Omega_0=\bigcap_{\ell\geq1}\Omega_\ell$
also has probability one, and all these convergences
hold simultaneously on $\Omega_0$.

Fix $\omega\in\Omega_0$. For each $\ell$, the
corresponding summand in
$d_{\rm w}(\cM_{s_k}(\omega),\mu_{\rm eq})$,
defined in \eqref{e:weak-metric}, tends to zero
as $k\to\infty$ and is bounded by $2^{-\ell}$.
Since $\sum_{\ell\geq1}2^{-\ell}<\infty$,
dominated convergence for the series gives
\[
d_{\rm w}(\cM_{s_k}(\omega),\mu_{\rm eq})
\longrightarrow0.
\]
Hence $\cM_{s_k}\to\mu_{\rm eq}$ weakly almost
surely. The inclusion of $f_1=1$ ensures, in
particular, convergence of the total masses.
No independence between different rows has been used.
\end{proof}

\begin{proof}[Proof of Theorem~\ref{T:fixed-basis-dimension}
\textup{(ii)}: Necessity]
Suppose that $\cM_{s_k}\to\mu_{\rm eq}$ weakly almost surely
in the prescribed sequence of bases, and assume that
the row vectors are pairwise independent.
Almost-sure convergence implies convergence
in probability, so part~\textup{(i)} gives
$\E|\xi|^2=1$.

We show that $\E|\xi|^{2+2/n}<\infty$.
Suppose, to the contrary, that
$\E|\xi|^{2+2/n}=\infty$, and set $V=|\xi|^2$ and
$p=1+1/n$. Then $\E V=1$, whereas $\E V^p=\infty$.
A single coefficient with squared modulus greater than $2N_k$
forces the total mass to exceed $2$. We will show that such a
coefficient occurs in infinitely many rows almost surely.

We first relate the moment of $V$ to its tail probabilities.
For any fixed $a>0$, Tonelli's theorem gives
\[
\sum_{k\geq1}k^n\P(V>ak^n)
=\E\left[\sum_{k\geq1}k^n\1_{\{V>ak^n\}}\right]
=\E\left[\sum_{1\leq k<(V/a)^{1/n}}k^n\right].
\]
For $t\geq2$, comparison with the integral of $x^n$ gives
\[
c_n t^{n+1}\leq\sum_{1\leq k<t}k^n\leq C_n t^{n+1},
\]
where $c_n,C_n>0$ depend only on $n$. Taking
$t=(V/a)^{1/n}$, we see that the inner sum is comparable to
$V^{1+1/n}$ for sufficiently large $V$. On bounded ranges
of $V$, both quantities are bounded, so their contributions
do not affect finiteness of the expectations. Consequently,
\begin{equation}\label{e:tail-moment-equivalence}
\sum_{k\geq1}k^n\P(V>ak^n)<\infty
\quad\Longleftrightarrow\quad \E V^p<\infty.
\end{equation}

By the dimension asymptotics, there are constants $c,C>0$
such that $ck^n\leq N_k\leq Ck^n$ for all sufficiently
large $k$. Hence
$N_k\P(V>2N_k)\geq ck^n\P(V>2Ck^n)$.
Since $\E V^p=\infty$, the equivalence
\eqref{e:tail-moment-equivalence} implies
\[
\sum_kN_k\P(V>2N_k)=\infty.
\]

Set $q_k=\P(V>2N_k)$. Let
$E_k=\{\max_{1\leq j\leq N_k}|\xi_{j,k}|^2>2N_k\}$.
Independence of the coordinates within each row gives
$\P(E_k)=1-(1-q_k)^{N_k}$. Moreover, since $\E V=1$,
\[
N_kq_k\leq\frac12\E\bigl(V\1_{\{V>2N_k\}}\bigr)
\longrightarrow0.
\]
The convergence follows from the integrability of $V$ and
$N_k\to\infty$. Using $(1-q_k)^{N_k}\leq e^{-N_kq_k}$
and $1-e^{-t}\geq t/2$ for $0\leq t\leq1$, we obtain
$\P(E_k)\geq N_kq_k/2$ for all sufficiently large $k$.
Consequently, $\sum_k\P(E_k)=\infty$.

Since the row vectors are pairwise independent, so are
the events $E_k$. As $\sum_k\P(E_k)=\infty$, the
pairwise-independent version of the second Borel--Cantelli
lemma \cite[Theorem~2.3.9]{Durrett2019} gives
\[
\P(E_k\text{ infinitely often})=1.
\]

On $E_k$, at least one coefficient satisfies
$|\xi_{j,k}|^2>2N_k$. By orthonormality,
\[
\cM_{s_k}(K)
=\frac1{N_k}\sum_{j=1}^{N_k}|\xi_{j,k}|^2>2.
\]
Thus $\cM_{s_k}(K)>2$ infinitely often almost surely,
contradicting the assumed convergence of the total masses to
$\mu_{\rm eq}(K)=1$. We conclude that
$\E|\xi|^{2+2/n}<\infty$.
\end{proof}

The same normalization argument yields the almost-sure
counterpart for self-normalized mass.

\begin{Corollary}
\label{C:self-normalized-as}
Suppose that
$0<\sigma^2:=\E|\xi|^2<\infty$, and
$\E|\xi|^{2+2/n}<\infty$.
Then $\widehat{\cM}_{s_k}\to\mu_{\rm eq}$ weakly
almost surely, without any independence assumption
between coefficient rows.
\end{Corollary}

\begin{proof}
Set $t_k=s_k/\sigma$. By the sufficiency
part of Theorem~\ref{T:fixed-basis-dimension}\textup{(ii)},
$\cM_{t_k}\to\mu_{\rm eq}$ weakly almost surely.

On the probability-one event where this convergence
holds, $\cM_{t_k}(K)\to1$, so the sections are
eventually nonzero. The normalization identity used
in the proof of Corollary~\ref{C:self-normalized-prob}
therefore gives
$\widehat{\cM}_{s_k}\to\mu_{\rm eq}$
weakly almost surely.
\end{proof}

\section{From mass equidistribution to zeros}
\label{S:MassZeros}

\subsection{A deterministic mass-to-zero implication}

In this section we assume that $X$ is K\"ahler. We retain
the notation of Section~\ref{S:Geometry}, in particular
$h=h_0e^{-\phi}$, $\theta=c_1(L,h_0)$, and the equilibrium
potential $P_K\phi$. We impose the additional condition
\begin{equation}\label{e:bulk-assumptions-section}
U:=\operatorname{Int}_X(\supp\mu_{\rm eq}),
\qquad \mu_{\rm eq}(U)=1.
\end{equation}
Thus the equilibrium measure gives full mass to the interior
of its support, which we call the open bulk. This condition
holds, for example, when $K=X$ and $h$ has smooth positive
curvature. It may fail when the equilibrium measure is
supported on a thin set, as Example~\ref{P:mass-not-zero}
shows.

A deterministic mass-to-zero implication for positive line
bundles appears in Shiffman--Zelditch
\cite[Lemma~1.4]{ShiffmanZelditch1999}, and a weighted
polynomial version was established by Bayraktar
\cite[Theorem~2.4]{Bayraktar2020}. We prove a corresponding
result for sections of big line bundles with continuous
metrics and Bernstein--Markov measures. 
In this setting, the equilibrium potential may be singular,
and an $L^1$ cluster limit of the potentials associated
with the normalized zero currents cannot be assumed
a priori to have minimal singularities. We first show
that each such limit agrees with $P_K\phi$ on the open
bulk and then extend this equality globally using
non-pluripolar Monge--Amp\`ere theory.

We begin with a consequence of the equilibrium contact
property.

\begin{Lemma}[Contact on the bulk]
\label{L:bulk-contact}
Under the above assumptions,
\begin{equation}\label{e:bulk-contact}
P_K\phi=\phi\qquad\text{on }U.
\end{equation}
In particular, $P_K\phi$ is continuous on $U$.
\end{Lemma}

\begin{proof}
By \cite[Proposition~2.10]{BermanBoucksom2010},
the measure $\mu_{\rm eq}$ is concentrated on $K$,
and $P_K\phi=\phi$ holds $\mu_{\rm eq}$-almost
everywhere. Since $K$ is closed,
$\supp\mu_{\rm eq}\subset K$, and hence $U\subset K$.

Let $Q_K\phi$ be the unregularized supremum
defining the equilibrium potential, so that $P_K\phi=(Q_K\phi)^*$.
 Every admissible potential is bounded above by
$\phi$ on $K$, and consequently $Q_K\phi\leq\phi$
on $U$. Since $U$ is open and $\phi$ is continuous,
taking upper limits within $U$ gives $P_K\phi\leq\phi$
there.

Suppose that $(P_K\phi)(x_0)<\phi(x_0)$ for some
$x_0\in U$. Upper semicontinuity of $P_K\phi$
and continuity of $\phi$ then give a neighborhood
$V\Subset U$ of $x_0$ and $\delta>0$ such that
$P_K\phi\leq\phi-\delta$ on $V$. But
$x_0\in\supp\mu_{\rm eq}$ implies that $\mu_{\rm eq}(V)>0$,
contradicting the almost-everywhere contact property.
\end{proof}

Next, we recall the pluripotential-theoretic facts needed
to identify the limiting potentials globally.
For $u\in\operatorname{PSH}(X,\theta)$, let
\begin{equation}\label{e:MA-theta-definition}
\operatorname{MA}_{\theta}(u)
:=\big\langle(\theta+dd^cu)^n\big\rangle
\end{equation}
denote its non-pluripolar Monge--Amp\`ere measure
\cite{BEGZ2010}. It agrees with the classical
Bedford--Taylor product wherever $u$ is locally bounded,
but is defined for arbitrary $\theta$-psh potentials
and does not charge pluripolar sets. 
Its total mass may be strictly smaller than $\vol(L)$,
but always satisfies
$\int_X\operatorname{MA}_{\theta}(u)\leq\vol(L)$
by \cite[Proposition~1.20]{BEGZ2010}. Following
\cite{BEGZ2010}, the \emph{full Monge--Amp\`ere mass
class} is given by
\begin{equation}\label{e:E-class}
\mathcal E(X,\theta)
:=\left\{u\in\operatorname{PSH}(X,\theta):
\int_X\operatorname{MA}_{\theta}(u)=\vol(L)\right\}.
\end{equation}
Potentials with minimal singularities have full
Monge--Amp\`ere mass. In particular,
$P_K\phi\in\mathcal E(X,\theta)$ and
$\operatorname{MA}_{\theta}(P_K\phi)
=\vol(L)\mu_{\rm eq}$.

We will use two properties of the non-pluripolar product.
The first is locality in the plurifine topology
\cite[Proposition~1.4]{BEGZ2010}. In particular,
if two $\theta$-psh potentials $u,v$ agree on
an open set $O\subset X$, then
\begin{equation}\label{e:locality-tool}
\1_O\operatorname{MA}_{\theta}(u)
=\1_O\operatorname{MA}_{\theta}(v).
\end{equation}

The second is the full-mass domination principle
of Darvas--Di Nezza--Lu
\cite[Proposition~2.4]{DDL2018}, based on an
argument of Dinew: if
$u\in\mathcal E(X,\theta)$ and
$v\in\operatorname{PSH}(X,\theta)$ satisfy
\begin{equation}\label{e:domination-tool}
\operatorname{MA}_{\theta}(u)(\{u<v\})=0,
\end{equation}
then $u\geq v$ on $X$. Note that no full-mass assumption
is required on $v$. In the proof below, locality
will show that a cluster potential has full mass,
allowing us to apply the domination principle
with that potential in the role of $u$.

\begin{Theorem}
\label{T:deterministic-mass-zero}
Assume~\eqref{e:bulk-assumptions-section}, and let
$s_k\in H^0(X,L^k)$ be a deterministic
sequence of nonzero sections, defined for all
sufficiently large $k$. If
\begin{equation}\label{e:deterministic-mass-convergence}
\cM_{s_k}\longrightarrow\mu_{\rm eq}
\end{equation}
weakly as finite measures on $X$, then
\begin{equation}\label{e:potential-convergence}
\frac1k\log|s_k|_{h_0^k}^2
\longrightarrow P_K\phi
\qquad\text{in }L^1(X).
\end{equation}
Equivalently, since $h=h_0e^{-\phi}$,
$k^{-1}\log|s_k|_{h^k}^2\to P_K\phi-\phi$
in $L^1(X)$. Moreover, we have
\begin{equation}\label{e:zero-convergence-deterministic}
\frac1k[Z_{s_k}]
\longrightarrow T_{\rm eq}=\theta+dd^cP_K\phi
\end{equation}
weakly as currents.
\end{Theorem}

\begin{proof}
Set $u_k=k^{-1}\log|s_k|_{h_0^k}^2$.
In a local holomorphic frame $e$ of $L$, write
$|e|_{h_0}^2=e^{-\rho}$ and
$s_k=f_ke^{\otimes k}$. Then
$u_k=k^{-1}\log|f_k|^2-\rho$, and the
Poincar\'e--Lelong formula gives
\begin{equation}\label{e:PL-zero}
\theta+dd^cu_k=\frac1k[Z_{s_k}].
\end{equation}
In particular, $u_k\in\operatorname{PSH}(X,\theta)$.
Since $h=h_0e^{-\phi}$, the mass measures satisfy
\begin{equation}\label{e:normalized-density-zero}
\cM_{s_k}=\frac1{N_k}e^{k(u_k-\phi)}\,d\mu.
\end{equation}
We show that $(u_k)$ is relatively compact in
$L^1(X)$ and that every cluster limit equals
$P_K\phi$.

Testing \eqref{e:deterministic-mass-convergence}
against the constant function $1$ gives
\begin{equation}\label{e:normalized-norm-zero}
\frac{\|s_k\|_k^2}{N_k}
=\cM_{s_k}(X)\longrightarrow1.
\end{equation}

Since $N_k\asymp k^n$, it follows that
$k^{-1}\log\|s_k\|_k^2\to0$. Fix $\varepsilon>0$, then the Bernstein--Markov inequality \eqref{e:BM} yields
\[
\sup_K(u_k-\phi)
\leq2\varepsilon+
\frac{2\log C_\varepsilon+\log\|s_k\|_k^2}{k}
=2\varepsilon+o(1).
\]
Since the right-hand side is independent of $x\in X$,
subtracting it from $u_k$ gives a $\theta$-psh
function bounded above by $\phi$ on $K$.
By the definition of the equilibrium envelope,
\begin{equation}\label{e:global-upper-zero}
u_k\leq P_K\phi+2\varepsilon+o(1)
\qquad\text{on }X.
\end{equation}
In particular, $(u_k)$ is uniformly bounded above.

We also show that $\sup_Xu_k$ remains bounded
from below. Indeed, since $\mu$ is a probability measure,
\[
\|s_k\|_k^2
=\int_K e^{k(u_k-\phi)}\,d\mu
\leq\exp\bigl(k(\sup_Xu_k-\inf_K\phi)\bigr).
\]
Taking logarithms and using that
$k^{-1}\log\|s_k\|_k^2\to0$, we obtain
$\sup_Xu_k\geq\inf_K\phi+o(1)$.
Now choose a K\"ahler form $\omega$ and $C>0$
such that $C\omega-\theta\geq0$. Then
$u_k-\sup_Xu_k\in\operatorname{PSH}(X,C\omega)$
has supremum zero. By
\cite[Proposition~2.7]{GuedjZeriahi2005},
such normalized potentials form a compact family
in $L^1(X)$. Since the normalizing constants
$\sup_Xu_k$ are bounded, $(u_k)$ is relatively
compact in $L^1(X)$.

Let $u_{k_j}\to G$ in $L^1(X)$ along a subsequence,
where $G\in\operatorname{PSH}(X,\theta)$ denotes
the upper-semicontinuous representative of the limit.
We first claim that
\begin{equation}\label{e:G-below-P}
G\leq P_K\phi\qquad\text{on }X.
\end{equation}
The argument is standard in pluripotential theory, but we include it for
completeness. Fix $\varepsilon>0$. Passing to a further
subsequence if necessary, we may assume
$u_{k_j}\to G$ almost everywhere. By
\eqref{e:global-upper-zero},
$G\leq P_K\phi+2\varepsilon$ almost everywhere.

To obtain the pointwise inequality, fix $x\in X$
and choose a smooth local potential $\rho$ of $\theta$
near $x$. The submean inequality for the
psh function $G+\rho$, together with the
almost-everywhere bound, gives, for sufficiently
small coordinate balls $B(x,r)$,
\[
G(x)+\rho(x)
\leq\frac1{\operatorname{Vol}(B(x,r))}
\int_{B(x,r)}(P_K\phi+\rho+2\varepsilon).
\]
By upper semicontinuity of $P_K\phi+\rho$,
the upper limit of the right-hand side as
$r\downarrow0$ is at most
$(P_K\phi)(x)+\rho(x)+2\varepsilon$.
Thus $G(x)\leq(P_K\phi)(x)+2\varepsilon$.
Letting $\varepsilon\downarrow0$
proves \eqref{e:G-below-P}.

We next prove equality on the open bulk.
Suppose that $G(x_0)<(P_K\phi)(x_0)$
for some $x_0\in U$. By
Lemma~\ref{L:bulk-contact},
$P_K\phi=\phi$ is continuous on $U$.
Upper semicontinuity of $G$ then allows
us to choose coordinate balls
$x_0\in B'\Subset B\Subset U$, a smooth
local potential $\rho$ of $\theta$ on a
neighborhood of $\overline B$, and
$\delta>0$ such that
\[
\sup_B(G+\rho)
\leq\inf_B(P_K\phi+\rho)-3\delta.
\]
The local Hartogs lemma, applied to the
$L^1$-convergent psh functions
$u_{k_j}+\rho$, gives
\[
\limsup_{j\to\infty}\sup_{B'}(u_{k_j}+\rho)
\leq\sup_B(G+\rho).
\]
It follows that, for all sufficiently large $j$,
$\sup_{B'}(u_{k_j}+\rho)
\leq\inf_B(P_K\phi+\rho)-2\delta$.
Consequently,
\[
u_{k_j}\leq P_K\phi-2\delta
\qquad\text{on }B'.
\]
Since $P_K\phi=\phi$ on $U$, the density
identity \eqref{e:normalized-density-zero}
yields
\begin{equation}\label{massdecay}
\begin{aligned}
\cM_{s_{k_j}}(B')
&=\frac1{N_{k_j}}
\int_{B'\cap K}e^{k_j(u_{k_j}-\phi)}\,d\mu \leq\frac{e^{-2k_j\delta}}{N_{k_j}}
\mu(B'\cap K)\longrightarrow0.
\end{aligned}
\end{equation}
On the other hand, $x_0\in\supp\mu_{\rm eq}$
implies $\mu_{\rm eq}(B')>0$. By the
Portmanteau theorem and
\eqref{e:deterministic-mass-convergence} we have
\[
\liminf_{j\to\infty}\cM_{s_{k_j}}(B')
\geq\mu_{\rm eq}(B')>0,
\]
which contradicts (\ref{massdecay}). We conclude that
\[
G=P_K\phi\qquad\text{on }U.
\]

It remains to extend this equality to $X$.
By locality \eqref{e:locality-tool}, we have
\[
\1_U\operatorname{MA}_{\theta}(G)
=\1_U\operatorname{MA}_{\theta}(P_K\phi).
\]
Since $\mu_{\rm eq}(U)=1$ and
$\operatorname{MA}_{\theta}(P_K\phi)
=\vol(L)\mu_{\rm eq}$, we obtain
$\int_U\operatorname{MA}_{\theta}(G)=\vol(L)$.
The universal upper bound on the total
non-pluripolar mass therefore implies that
$G\in\mathcal E(X,\theta)$ and
\[
\operatorname{MA}_{\theta}(G)(X\setminus U)=0.
\]
Since $G=P_K\phi$ on $U$, the set
$\{G<P_K\phi\}$ is contained in $X\setminus U$.
Hence
$\operatorname{MA}_{\theta}(G)(\{G<P_K\phi\})=0$.
The full-mass domination principle
\eqref{e:domination-tool}, applied with
$u=G$ and $v=P_K\phi$, gives
$G\geq P_K\phi$ on $X$. Together with
\eqref{e:G-below-P}, this proves
$G=P_K\phi$ everywhere.

Every $L^1$ cluster limit of $(u_k)$ is
therefore $P_K\phi$. Relative compactness
implies that the full sequence converges
to $P_K\phi$ in $L^1(X)$, proving
\eqref{e:potential-convergence}. Since
$h=h_0e^{-\phi}$, the equivalent convergence
of logarithmic potentials with respect to
$h$ follows by subtracting $\phi$.
Finally, \eqref{e:PL-zero} and continuity
of $dd^c$ on distributions give the
zero-current convergence
\eqref{e:zero-convergence-deterministic}.
\end{proof}

\begin{Remark}
The zero-current conclusion of
Theorem~\ref{T:deterministic-mass-zero}
also holds if the nonzero sections satisfy
$\widehat{\cM}_{s_k}\to\mu_{\rm eq}$
weakly. Indeed, the positive definiteness
of the $L^2$ inner product ensures that
$\|s_k\|_k>0$. Setting
$t_k=\sqrt{N_k}s_k/\|s_k\|_k$, we have
$\cM_{t_k}=\widehat{\cM}_{s_k}$ and
$[Z_{t_k}]=[Z_{s_k}]$. Applying the
theorem to $(t_k)$ gives
$k^{-1}[Z_{s_k}]\to T_{\rm eq}$ weakly,
as well as
\[
\frac1k\log|t_k|_{h_0^k}^2
\longrightarrow P_K\phi
\qquad\text{in }L^1(X).
\]
The corresponding convergence of logarithmic
potentials need not hold for the original
sections. Indeed, multiplying $s_k$ by
$e^{ck}$, with $c\in\R$, leaves its
self-normalized mass unchanged but adds
$2c$ to $k^{-1}\log|s_k|_{h_0^k}^2$.
\end{Remark}

\begin{Remark}
\label{R:non-kahler}
We state Theorem~\ref{T:deterministic-mass-zero}
in the K\"ahler setting for simplicity to apply the domination principle
directly on $X$. However, the conclusion remains valid without
this assumption.

Indeed, since $L$ is big, $X$ is Moishezon and admits
a modification $\pi:Y\to X$ with $Y$ smooth and
projective. 
The compactness argument applies to
$(u_k\circ\pi)$ on $Y$, and relative compactness
descends to $L^1(X)$ since $\pi$ has degree one
and $\pi^*dV_X\leq C\,dV_Y$ for some $C>0$
and smooth positive volume forms $dV_X,dV_Y$.
As before, every cluster
limit $G$ satisfies $G\leq P_K\phi$ on $X$ and
$G=P_K\phi$ on $U$.

Write $P=P_K\phi$. By locality and the modification
formula \cite[Remark~1.7]{BEGZ2010}, for
$v\in\{G,P\}$,
\[
\int_{\pi^{-1}(U)}
\operatorname{MA}_{\pi^*\theta}(v\circ\pi)
=\int_U\operatorname{MA}_\theta(P)
=\vol(L)=\vol(\pi^*L).
\]
Thus $G\circ\pi$ has full Monge--Amp\`ere mass,
concentrated on $\pi^{-1}(U)$, where
$G\circ\pi=P\circ\pi$. The domination principle
on the K\"ahler manifold $Y$ gives
$G\circ\pi\geq P\circ\pi$, and hence $G=P$ on $X$.
\end{Remark}

\subsection{Random zero currents and the bulk obstruction}

We now apply the deterministic result to random sections.
For a nonzero section $s_k$, the Poincar\'e--Lelong
formula shows that $k^{-1}[Z_{s_k}]$ is a positive
closed $(1,1)$-current representing $c_1(L)$. We therefore consider
the space of positive closed $(1,1)$-currents on $X$ in this
cohomology class,
\[
\mathscr T_L
:=\{T\geq0:\ dT=0,\ [T]=c_1(L)\},
\]
endowed with the weak topology of currents. Recall that
$T_j\to T$ weakly if
$\langle T_j,\psi\rangle\to\langle T,\psi\rangle$
for every smooth $(n-1,n-1)$-form $\psi$ on $X$.

This topology is metrizable on $\mathscr T_L$. Indeed,
for any K\"ahler form $\omega$ and $T\in\mathscr T_L$,
\[
\int_X T\wedge\omega^{n-1}
=\int_X\theta\wedge\omega^{n-1}.
\]
Thus the currents in $\mathscr T_L$ have uniformly bounded
mass. By positivity, weak convergence is determined by a
countable dense family of smooth test forms. We use the
resulting notions of convergence in probability and almost
surely for random currents taking values in $\mathscr T_L$.

On the event $\{s_k=0\}$, we set
$k^{-1}[Z_{s_k}]=T_{\rm eq}$. This convention does not
affect the limits, since the sections are eventually
nonzero almost surely, as shown in the proof below.

\begin{Corollary}\label{C:random-zero}
Assume~\eqref{e:bulk-assumptions-section} and
$0<\E|\xi|^2<\infty$. Then
\[
\frac1k[Z_{s_k}]\longrightarrow T_{\rm eq}
\]
in probability in the weak sense of currents. If, in addition,
$\E|\xi|^{2+2/n}<\infty$, the convergence is almost sure.
\end{Corollary}

\begin{proof}
Let $\sigma^2=\E|\xi|^2$ and set $t_k=s_k/\sigma$.
The coefficients of $t_k$ are centered and have unit
second moment, while $[Z_{t_k}]=[Z_{s_k}]$ whenever
$s_k\ne0$. We may therefore assume $\E|\xi|^2=1$.

Set $q=\P(\xi=0)<1$. Independence of the coefficients
within each row and linear independence of the basis
sections give $\P(s_k=0)=q^{N_k}$. Since $N_k\asymp k^n$,
\[
\sum_{k\geq1}\P(s_k=0)
=\sum_{k\geq1}q^{N_k}<\infty.
\]
By the first Borel--Cantelli lemma, the sections are
nonzero for all sufficiently large $k$ almost surely,
regardless of the dependence between rows.

Suppose first that $\E|\xi|^{2+2/n}<\infty$.
The sufficient direction of
Theorem~\ref{T:fixed-basis-dimension}\textup{(ii)}
gives $\cM_{s_k}\to\mu_{\rm eq}$ weakly almost surely.
For every outcome in the probability-one event where
this convergence holds and the sections are eventually
nonzero, we may apply
Theorem~\ref{T:deterministic-mass-zero} to the resulting
deterministic sequence of sections. It follows that
\[
\frac1k[Z_{s_k}]\longrightarrow T_{\rm eq}
\]
almost surely in the weak sense of currents.

To prove convergence in probability, let $(k_j)$ be any
subsequence. By
Theorem~\ref{T:fixed-basis-dimension}\textup{(i)},
$\cM_{s_{k_j}}\to\mu_{\rm eq}$ weakly in probability.
Since the space of finite measures on the compact set
$K$ is metrizable in the weak topology, there is a
further subsequence $(k_{j_\ell})$ along which the
mass measures converge weakly almost surely.
As these measures are supported on $K$, the convergence
also holds as measures on $X$.

On the probability-one event where the mass measures
converge and the sections are eventually nonzero,
the argument of Theorem~\ref{T:deterministic-mass-zero}
applies along this subsequence and yields
\[
\frac1{k_{j_\ell}}[Z_{s_{k_{j_\ell}}}]
\longrightarrow T_{\rm eq}
\]
weakly as currents. Thus every subsequence of the
random zero currents admits a further subsequence
converging almost surely to $T_{\rm eq}$.
Since these currents take values in the metrizable
space $\mathscr T_L$, the subsequence criterion for
convergence in probability completes the proof.
\end{proof}

Almost-sure zero equidistribution for big line bundles with
continuous metrics was established by
Bayraktar--Coman--Marinescu--Nguy\^en
\cite[Theorem~1.2]{BCMN2024} under a logarithmic-moment
condition on linear functionals of the coefficient vectors,
together with a suitable summability assumption. Their
logarithmic-moment condition need not hold for discrete
coefficient laws. For example, if $\xi_1$ and $\xi_2$
are independent Rademacher variables, then
$\P((\xi_1-\xi_2)/\sqrt{2}=0)=1/2$, so the corresponding
logarithmic moment is infinite.

Corollary~\ref{C:random-zero} applies to this ensemble under the
bulk hypothesis. More generally, it gives convergence in probability
under a finite positive second moment and almost-sure convergence
under $\E|\xi|^{2+2/n}<\infty$, without independence between
rows. The additional bulk hypothesis makes this result complementary
to~\cite{BCMN2024}, which treats more general complex spaces.
Our moment assumptions are sufficient for zero equidistribution;
the total-mass argument proving their necessity for
$\cM_{s_k}$ does not give a converse for zeros.

Next we provide an example, adapted from
\cite[Example~2.5]{Bayraktar2020}, which shows that mass equidistribution
alone need not determine the limiting zero current.

\begin{Example}[Mass convergence without zero equidistribution]
\label{P:mass-not-zero}
There exist a compact K\"ahler manifold $X$, a big line bundle
$L\to X$ with a smooth positively curved Hermitian metric $h$, a
Bernstein--Markov probability measure $\mu$ with non-pluripolar
support, and nonzero sections $s_k\in H^0(X,L^k)$ such that
\[
\cM_{s_k}=\mu_{\rm eq}\qquad\text{for every }k,
\]
but $k^{-1}[Z_{s_k}]$ does not converge to $T_{\rm eq}$.
\end{Example}

\begin{proof}
Take $X=\mathbf P^1$, $L=\mathcal O(1)$, and
$K=S^1\subset\C\subset\mathbf P^1$. Let $\mu$ be normalized
arclength on $S^1$, so that $d\mu(e^{i\vartheta})=d\vartheta/(2\pi)$.
Let $h_{\rm FS}$ be the Fubini--Study metric and set $h=2h_{\rm FS}$.
In homogeneous coordinates $[Z_0:Z_1]$, consider
the affine chart $Z_0\ne0$ with $z=Z_1/Z_0$.
The standard frame $e_0$ of $\mathcal O(1)$
satisfies
$|e_0|_h^2=2/(1+|z|^2)=1$ on $S^1$.

The Bernstein--Markov property follows from Parseval's identity
and the Cauchy--Schwarz inequality. Indeed, if
$s=p\,e_0^{\otimes k}$, where
$p(z)=\sum_{j=0}^ka_jz^j$, then
\[
\|s\|_k^2=\sum_{j=0}^k|a_j|^2,
\qquad
\sup_{S^1}|s|_{h^k}
\leq\sqrt{k+1}\,\|s\|_k.
\]
The factor $\sqrt{k+1}$ is subexponential in $k$, as required.

With $h_0=h_{\rm FS}$, the weight is $\phi=-\log2$. The
equilibrium potential in the affine chart is given by
\[
P_K\phi=2\log^+|z|-\log(1+|z|^2).
\]

Indeed, the function $2\log^+|z|$ is subharmonic,
vanishes on $S^1$, and has the prescribed logarithmic
growth at infinity. Conversely, the maximum principle
gives the upper bound for any admissible local weight
on the unit disk. Applying the same principle in the
chart at infinity, after subtracting $2\log|z|$,
gives the corresponding bound outside the disk.
Thus, in our normalization of $dd^c$,
\[
T_{\rm eq}
=\theta+dd^cP_K\phi
=dd^c(2\log^+|z|)
=\frac{d\vartheta}{2\pi}
\]
on the affine chart.
 In the coordinate $w=1/z$ near infinity,
$P_K\phi=-\log(1+|w|^2)$, so $T_{\rm eq}$ has no mass at
infinity. Consequently, $T_{\rm eq}=\mu_{\rm eq}=\mu$
on $\mathbf P^1$.

Since $N_k=k+1$, take $s_k=\sqrt{k+1}\,Z_1^k$. On $S^1$,
$|Z_1^k|_{h^k}^2=1$, and therefore
\[
\cM_{s_k}=\frac1{k+1}|s_k|_{h^k}^2\,d\mu
=\mu=\mu_{\rm eq}\qquad\text{for every }k.
\]
On the other hand, $[Z_{s_k}]=k[0]$, so
\[
\frac1k[Z_{s_k}]=[0]\ne T_{\rm eq}
\qquad\text{for every }k.
\]
Here $\supp\mu_{\rm eq}=S^1$ has empty interior in $\mathbf P^1$,
so~\eqref{e:bulk-assumptions-section} fails. The obstruction thus
occurs even for a smooth positive metric on an ample line bundle.
This does not show that the bulk hypothesis is necessary.
\end{proof}

\section{Fluctuations of mass statistics}\label{S:CLT}

We now study the fluctuations of the mass statistics $\cM_{s_k}(f)$ for
$f\in \mathcal C^0(K,\R)$. Throughout this section, the coefficient law satisfies
\begin{equation}\label{e:CLT-moments}
\E\xi=0,\qquad \E|\xi|^2=1,\qquad \nu_4:=\E|\xi|^4<\infty,
\end{equation}
and we put
\begin{equation}\label{e:beta-kappa}
\beta:=\E\xi^2,\qquad \kappa:=\nu_4-2-|\beta|^2.
\end{equation}
We impose neither any symmetry assumption nor the condition $\beta=0$.
For $f\in \mathcal C^0(K,\R)$, recall that the Toeplitz matrix $A_k(f)$ is Hermitian, and
\begin{equation}\label{e:CLT-quadratic-representation}
\cM_{s_k}(f)=\frac{\xi_k^*A_k(f)\xi_k}{N_k},\qquad
\E\cM_{s_k}(f)=\frac{\Tr A_k(f)}{N_k}.
\end{equation}

In particular,
\[
\sqrt{N_k}\bigl(\cM_{s_k}(f)-\E\cM_{s_k}(f)\bigr)
=\frac{\xi_k^*A_k(f)\xi_k-\Tr A_k(f)}{\sqrt{N_k}}.
\]
The fluctuations are thus governed by centered Hermitian
quadratic forms in independent complex random variables.
We study them using the central limit theorem for sesquilinear
forms of Bai--Yao~\cite{BaiYao2008,BaiYao2021}; see also
Bhansali--Giraitis--Kokoszka~\cite{BGK2007} for related
results in the real setting. On the geometric side, we use
the following Toeplitz estimates, established in
Section~\ref{S:Geometry}:
\begin{equation}\label{e:CLT-Toeplitz-input}
\|A_k(f)\|_{\rm op}\leq\|f\|_\infty,\qquad
\frac{\Tr\bigl(A_k(f)A_k(g)\bigr)}{N_k}
\longrightarrow\mu_{\rm eq}(fg),\quad f,g\in \mathcal C^0(K,\R).
\end{equation}

These estimates alone, however, do not determine the
asymptotic variance of $\cM_{s_k}(f)$ for general coefficient
laws. As the exact covariance formula in
Lemma~\ref{P:mass-covariance} shows, the variance also
involves the sum of squared diagonal entries and the
transpose trace. These quantities depend on the prescribed
orthonormal bases, and their normalized limits need not exist.

We therefore begin with the exact covariance and prove the
CLT by a subsequence argument, using the finite-degree
variance as normalization without requiring convergence of
the basis-dependent terms. For self-normalized mass, we
express the statistic as a ratio of quadratic forms and
isolate a centered quadratic form that governs its
fluctuations. A separate estimate then permits centering
at the statistic's own expectation. Finally, the
projective-line example illustrates the effect of the
choice of basis, while Proposition~\ref{T:Haar-basis-CLT}
identifies the limiting variances for Haar-typical bases.

\subsection{Covariance and variance bounds}

For Hermitian $d\times d$ matrices $A=(a_{ij})$ and $B=(b_{ij})$, set
\begin{equation}\label{e:covariance-form}
\mathcal V_\xi(A,B)
:=\Tr(AB)+|\beta|^2\Tr(AB^T)+\kappa\sum_{j=1}^d a_{jj}b_{jj}.
\end{equation}
By Hermitian symmetry, this is a real symmetric bilinear
form on the real vector space of Hermitian matrices.

\begin{Lemma}[Exact covariance]\label{P:mass-covariance}
If $z$ has independent coordinates distributed as $\xi$, then
\[
\Cov(z^*Az,z^*Bz)=\mathcal V_\xi(A,B).
\]
Consequently, for $f,g\in \mathcal C^0(K,\R)$,
\begin{equation}\label{e:mass-covariance-general}
N_k^2\Cov\bigl(\cM_{s_k}(f),\cM_{s_k}(g)\bigr)
=\mathcal V_\xi\bigl(A_k(f),A_k(g)\bigr).
\end{equation}
\end{Lemma}

\begin{proof}
Expanding the quadratic forms gives
\[
\E\bigl[(z^*Az)(z^*Bz)\bigr]
=\sum_{i,j,p,q=1}^d a_{ij}b_{pq}
\E(\ov\xi_i\xi_j\ov\xi_p\xi_q).
\]
By independence and centering, an index which occurs exactly once
makes the corresponding expectation vanish. The remaining terms come
from the three pairings of the four indices, together with a correction
when they are all equal. Explicitly,
\begin{align*}
\E(\ov\xi_i\xi_j\ov\xi_p\xi_q)
={}&\delta_{ij}\delta_{pq}+\delta_{iq}\delta_{jp}
   +|\beta|^2\delta_{ip}\delta_{jq}
   +\kappa\delta_{ij}\delta_{jp}\delta_{pq}.
\end{align*}
The pairings $i=j$, $p=q$ and $i=q$, $j=p$ each contribute $1$,
whereas $i=p$, $j=q$ contributes
$\E\ov\xi^{\,2}\,\E\xi^2=|\beta|^2$. When all indices agree,
the sum of these three contributions is $2+|\beta|^2$, which the
last term corrects to $\nu_4$.

After multiplication by $a_{ij}b_{pq}$ and summation, the four terms
become, in the same order,
\[
\Tr A\,\Tr B,\qquad \Tr(AB),\qquad
|\beta|^2\Tr(AB^T),\qquad \kappa\sum_j a_{jj}b_{jj}.
\]
Since $\E(z^*Az)=\Tr A$ and $\E(z^*Bz)=\Tr B$, subtracting the
product of the expectations proves the first assertion, and the
quadratic-form representation~\eqref{e:CLT-quadratic-representation}
gives the second assertion.
\end{proof}

Taking $A=B$ and using $A^*=A$ in~\eqref{e:covariance-form}
yields
\begin{equation}\label{e:mass-variance-split}
\mathcal V_\xi(A,A)
=2\sum_{i<j}\bigl(|a_{ij}|^2+|\beta|^2\operatorname{Re}(a_{ij}^2)\bigr)
 +(\nu_4-1)\sum_j a_{jj}^2.
\end{equation}
In particular,
$\Tr(AA^T)=\sum_j a_{jj}^2+2\sum_{i<j}\operatorname{Re}(a_{ij}^2)$
is real, but need not be nonnegative. Thus the diagonal and off-diagonal
parts of the matrix make separate contributions to the variance.

\begin{Lemma}[Variance comparison]\label{L:mass-variance-comparison}
There is a constant $C_\xi>0$, depending only on the law of $\xi$,
such that, for every Hermitian matrix $A$,
\begin{equation}\label{e:mass-variance-upper}
0\leq\mathcal V_\xi(A,A)\leq C_\xi\Tr(A^2).
\end{equation}
If
\begin{equation}\label{e:CLT-nondegeneracy}
|\beta|<1,\qquad \nu_4>1,
\end{equation}
then there is also a constant $c_\xi>0$ such that
\begin{equation}\label{e:mass-variance-comparison-general}
c_\xi\Tr(A^2)\leq\mathcal V_\xi(A,A)\leq C_\xi\Tr(A^2).
\end{equation}
One may take
$C_\xi=\max\{1+|\beta|^2,\nu_4-1\}$ and
$c_\xi=\min\{1-|\beta|^2,\nu_4-1\}$.
Conversely, the existence of a positive lower constant valid for
all Hermitian matrices implies~\eqref{e:CLT-nondegeneracy}.
\end{Lemma}

\begin{proof}
For every $a\in\C$, the elementary estimate
\[
(1-|\beta|^2)|a|^2
\leq |a|^2+|\beta|^2\operatorname{Re}(a^2)
\leq(1+|\beta|^2)|a|^2
\]
follows from $|\operatorname{Re}(a^2)|\leq|a|^2$. Since
$|\beta|\leq\E|\xi|^2=1$ and
$\nu_4-1=\Var(|\xi|^2)\geq0$, applying these bounds to
\eqref{e:mass-variance-split} and using
\[
\Tr(A^2)=2\sum_{i<j}|a_{ij}|^2+\sum_j a_{jj}^2
\]
gives~\eqref{e:mass-variance-upper}. Under
\eqref{e:CLT-nondegeneracy}, the lower coefficients are both strictly
positive, and~\eqref{e:mass-variance-comparison-general} follows.

For the converse, a nonzero diagonal matrix has zero variance when
$\nu_4=1$. If $|\beta|=1$,  the Hermitian matrix
$\left(\begin{smallmatrix}0&i\\-i&0\end{smallmatrix}\right)$
has zero variance by~\eqref{e:mass-variance-split}, although its
Hilbert--Schmidt square is $2$. Either case precludes a positive
uniform lower bound.
\end{proof}

\subsection{Central limit theorem (CLT)}

The covariance formula involves three matrix quantities: the
Hilbert--Schmidt square, the sum of squared diagonal entries, and
the transpose trace. When their normalized limits exist, the
following quadratic-form theorem gives the Gaussian limit.

All Gaussian limits in what follows are real. We write
$\mathcal N(0,v)$ for the centered real normal law of variance $v\geq0$,
with $\mathcal N(0,0)$ understood as the point mass at zero.
We use the specialization of the central limit theorem of Bai and
Yao~\cite[Theorem~7.1]{BaiYao2008} to a single Hermitian quadratic
form, with the uniform operator-norm bound specified in the
correction~\cite{BaiYao2021}.

\begin{Theorem}[Bai--Yao]\label{T:Bai-Yao-mass-form}
Let $B_k=(b_{ij}^{(k)})$ be deterministic Hermitian matrices of sizes
$d_k\to\infty$, and suppose that $\sup_k\|B_k\|_{\rm op}<\infty$ and
\begin{equation}\label{e:BY-profiles}
\frac{\Tr(B_k^2)}{d_k}\to\vartheta,\qquad
\frac1{d_k}\sum_j(b_{jj}^{(k)})^2\to\omega,\qquad
\frac{\Tr(B_kB_k^T)}{d_k}\to\tau.
\end{equation}
If $z_k$ has independent coordinates with the law of $\xi$, then
\begin{equation}\label{e:BY-conclusion}
\frac{z_k^*B_kz_k-\Tr B_k}{\sqrt{d_k}}
\longrightarrow\mathcal N(0,\vartheta+|\beta|^2\tau+\kappa\omega).
\end{equation}
The convergence is in distribution, and the limiting variance may be zero.
\end{Theorem}

Indeed, in the notation of Bai--Yao, taking a single form with $X=Y$
and $\E|\xi|^2=1$ gives the three variance contributions
$\omega(\nu_4-1)$, $\vartheta-\omega$, and
$|\beta|^2(\tau-\omega)$. Their sum yields the variance
in~\eqref{e:BY-conclusion}. No independence of the coefficient vectors
for different values of $k$ is assumed or needed.

For the Toeplitz matrices, the geometric estimates
give convergence of the first profile, while the other
two remain bounded. The proof of
Theorem~\ref{T:intro-exact-variance-CLT} therefore
proceeds by extracting convergent profiles along
subsequences and normalizing by the exact variance.

\begin{proof}[Proof of Theorem~\ref{T:intro-exact-variance-CLT}]
Fix $f\in\mathcal C^0(K,\R)$, and write
$A_k=A_k(f)=(a_{ij}^{(k)})$.
By~\eqref{e:CLT-Toeplitz-input}, we have 
\begin{equation}\label{e:CLT-main-trace}
\frac{\Tr(A_k^2)}{N_k}\longrightarrow\mu_{\rm eq}(f^2),
\qquad
\sup_k\|A_k\|_{\rm op}\leq\|f\|_\infty.
\end{equation}
The other two profiles in~\eqref{e:BY-profiles}
need not converge, but they are uniformly bounded.
Indeed, since $A_k$ is Hermitian,
\[
0\leq\frac1{N_k}\sum_j(a_{jj}^{(k)})^2
\leq\frac{\Tr(A_k^2)}{N_k},
\qquad
\left|\frac{\Tr(A_kA_k^T)}{N_k}\right|
\leq\frac{\Tr(A_k^2)}{N_k}.
\]

The first bound follows from
$\Tr(A_k^2)=\sum_{i,j}|a_{ij}^{(k)}|^2$,
while the second follows from
$\Tr(A_kA_k^T)=\sum_{i,j}(a_{ij}^{(k)})^2$
and the triangle inequality. Moreover, since
\[
\frac{\Tr(A_k^2)}{N_k}
\leq\|A_k\|_{\rm op}^2
\leq\|f\|_\infty^2,
\]
both profiles are uniformly bounded.
Let $(k_m)$ be any subsequence. By the preceding
bounds, we may pass to a further subsequence,
still denoted by $(k_m)$, such that
\[
\frac1{N_{k_m}}\sum_j(a_{jj}^{(k_m)})^2
\longrightarrow\omega,
\qquad
\frac{\Tr(A_{k_m}A_{k_m}^T)}{N_{k_m}}
\longrightarrow\tau
\]
for some $\omega,\tau\in\R$.
The first profile converges along the same
subsequence to $\mu_{\rm eq}(f^2)$
by~\eqref{e:CLT-main-trace}.
Theorem~\ref{T:Bai-Yao-mass-form} therefore gives
\begin{equation}\label{e:CLT-subsequence-limit}
\sqrt{N_{k_m}}\bigl(
\cM_{s_{k_m}}(f)-\E\cM_{s_{k_m}}(f)
\bigr)
\longrightarrow\mathcal N(0,\sigma^2)
\end{equation}
in distribution, where $\sigma^2
=\mu_{\rm eq}(f^2)+|\beta|^2\tau+\kappa\omega.
$
On the other hand, the exact covariance formula yields
\[
\begin{aligned}
N_k\Var\bigl(\cM_{s_k}(f)\bigr)
={}&\frac{\Tr(A_k^2)}{N_k}
+\frac{|\beta|^2}{N_k}\Tr(A_kA_k^T)+\frac{\kappa}{N_k}
\sum_j(a_{jj}^{(k)})^2.
\end{aligned}
\]
Hence
\begin{equation}\label{e:CLT-subsequence-variance}
N_{k_m}\Var\bigl(\cM_{s_{k_m}}(f)\bigr)
\longrightarrow\sigma^2.
\end{equation}
By~\eqref{e:CLT-variance-lower-intro},
we have $\sigma^2>0$. Combining
\eqref{e:CLT-subsequence-limit} and
\eqref{e:CLT-subsequence-variance}
with Slutsky's theorem, we obtain
\[
\frac{\cM_{s_{k_m}}(f)-\E\cM_{s_{k_m}}(f)}
{\sqrt{\Var(\cM_{s_{k_m}}(f))}}
\longrightarrow\mathcal N(0,1)
\]
in distribution. Since every subsequence
admits a further subsequence with this same
limit, convergence holds along the full sequence:
\begin{equation}\label{e:mass-exact-variance-CLT}
\frac{\cM_{s_k}(f)-\E\cM_{s_k}(f)}
{\sqrt{\Var(\cM_{s_k}(f))}}
\longrightarrow\mathcal N(0,1).
\end{equation}

Finally, suppose that~\eqref{e:CLT-nondegeneracy-intro}
holds. By Lemma~\ref{L:mass-variance-comparison},
\[
N_k\Var\bigl(\cM_{s_k}(f)\bigr)
=\frac{\mathcal V_\xi(A_k,A_k)}{N_k}
\geq c_\xi\frac{\Tr(A_k^2)}{N_k}.
\]
Taking lower limits and using~\eqref{e:CLT-main-trace},
we obtain
\[
\liminf_{k\to\infty}
N_k\Var\bigl(\cM_{s_k}(f)\bigr)
\geq c_\xi\mu_{\rm eq}(f^2)>0.
\]
Thus~\eqref{e:CLT-variance-lower-intro}
holds for every prescribed sequence of
orthonormal bases, completing the proof.
\end{proof}

\begin{Remark}[Nondegeneracy and normalization]\label{R:degenerate-laws}
The coefficient conditions~\eqref{e:CLT-nondegeneracy} are sufficient
for the variance lower bound, but not necessary. If $f=1$, then
\[
\cM_{s_k}(1)=\frac1{N_k}\sum_j|\xi_{j,k}|^2,
\qquad N_k\Var\bigl(\cM_{s_k}(1)\bigr)=\nu_4-1,
\]
so a nondegenerate CLT holds for any $\nu_4>1$, including real
coefficient laws with $|\beta|=1$. On the other hand, when $|\xi|=1$
almost surely, a diagonal quadratic form is deterministic, although
an off-diagonal one need not be. This distinction will be visible
in the projective-line example below.

If the diagonal-square and transpose-trace profiles do converge for
a particular basis sequence, Bai--Yao gives a CLT on the fixed
$\sqrt{N_k}$ scale with the corresponding explicit limiting variance.
The point of~\eqref{e:mass-exact-variance-CLT} is that these
additional profile limits are not required after normalization by
the exact variance.
\end{Remark}

\subsection{Self-normalized fluctuations}

We now consider the self-normalized mass statistic
$\widehat{\cM}_{s_k}(f)$, which, on $\{s_k\ne0\}$,
is the ratio of two quadratic forms:
\[
\widehat{\cM}_{s_k}(f)
=\frac{\xi_k^*A_k(f)\xi_k}{\xi_k^*\xi_k}.
\]
The Bai--Yao theorem applies to quadratic forms, but not
directly to their ratios. To study the fluctuations of this
ratio, we subtract the quotient of the expectations of its
numerator and denominator. Since these expectations are
$\Tr A_k(f)$ and $N_k$, respectively, their quotient is
$\mu_k^{\mathrm{Berg},\Delta}(f)$. This leads us to consider the
trace-free Hermitian matrix
\begin{equation}\label{e:finite-k-centering}
F_k(f):=A_k(f)-\mu_k^{\mathrm{Berg},\Delta}(f)I_{N_k},
\end{equation}
where $I_{N_k}$ denotes the identity matrix.
By construction, $\Tr F_k(f)=0$. Moreover, since
$\mu_k^{\mathrm{Berg},\Delta}$ is a probability measure,
the operator-norm bound~\eqref{e:CLT-Toeplitz-input}
gives $\|F_k(f)\|_{\rm op}\leq2\|f\|_\infty$.

Expanding the product and taking traces, we obtain
\[
\frac{\Tr(F_k(f)F_k(g))}{N_k}
=\frac{\Tr(A_k(f)A_k(g))}{N_k}
-\mu_k^{\mathrm{Berg},\Delta}(f)
\mu_k^{\mathrm{Berg},\Delta}(g).
\]
The diagonal Bergman measure limit and
\eqref{e:CLT-Toeplitz-input} therefore give
\begin{equation}\label{e:centered-trace-limit}
\frac{\Tr(F_k(f)F_k(g))}{N_k}
\longrightarrow
\mu_{\rm eq}(fg)-\mu_{\rm eq}(f)\mu_{\rm eq}(g)
\end{equation}
for $f,g\in \mathcal C^0(K,\R)$.

Recall that $r_k=\|s_k\|_k^2/N_k
=\xi_k^*\xi_k/N_k$. By the definition of $F_k(f)$,
whenever $s_k\ne0$, we have
\begin{equation}\label{e:ratio-identity}
\sqrt{N_k}\bigl(\widehat{\cM}_{s_k}(f)
-\mu_k^{\mathrm{Berg},\Delta}(f)\bigr)
=\frac{\xi_k^*F_k(f)\xi_k}{\sqrt{N_k}\,r_k}.
\end{equation}
Since $\Tr F_k(f)=0$, the normalized quadratic form
$N_k^{-1/2}\xi_k^*F_k(f)\xi_k$ is centered.
By Lemma~\ref{P:mass-covariance}, its variance is
\begin{equation}\label{e:self-normalized-variance}
v_k (f):=\frac1{N_k}
\mathcal V_\xi\bigl(F_k(f),F_k(f)\bigr).
\end{equation}
This provides the variance scale for the self-normalized
fluctuations, although it need not be the
exact variance of the ratio itself.

\begin{Theorem}[Self-normalized CLT]
\label{T:self-normalized-CLT}
Assume the fourth-moment hypotheses~\eqref{e:CLT-moments}, and
fix $f\in \mathcal C^0(K,\R)$.
If $\liminf_{k\to\infty}v_k (f)>0$, then
\begin{equation}\label{e:self-normalized-CLT}
\frac{\sqrt{N_k}\bigl(\widehat{\cM}_{s_k}(f)
-\E\widehat{\cM}_{s_k}(f)\bigr)}
{\sqrt{v_k (f)}}
\longrightarrow\mathcal N(0,1)
\end{equation}
in distribution. If $v_k (f)\to v\geq0$, then
\begin{equation}\label{e:self-normalized-CLT-degenerate}
\sqrt{N_k}\bigl(\widehat{\cM}_{s_k}(f)
-\E\widehat{\cM}_{s_k}(f)\bigr)
\longrightarrow\mathcal N(0,v)
\end{equation}
in distribution. 
In either assertion, the expectation may be replaced
by $\mu_k^{\mathrm{Berg},\Delta}(f)$.
Moreover, the positive variance lower bound holds
under~\eqref{e:CLT-nondegeneracy} whenever
$\mu_{\rm eq}(f^2)>\mu_{\rm eq}(f)^2$.
\end{Theorem}

\begin{proof}
Consider the centered quadratic forms
\[
W_k:=\frac{\xi_k^*F_k(f)\xi_k}{\sqrt{N_k}}.
\]
By Lemma~\ref{P:mass-covariance},
$\E W_k=0$ and $\E|W_k|^2=v_k(f)$.

We first establish the corresponding central limit
assertions for $W_k$. The argument follows the proof of
Theorem~\ref{T:intro-exact-variance-CLT}, with
$F_k(f)$ in place of $A_k(f)$.
Indeed, these matrices have uniformly bounded operator
norms by~\eqref{e:CLT-Toeplitz-input} and the fact that
$\mu_k^{\mathrm{Berg},\Delta}$ is a probability measure,
while~\eqref{e:centered-trace-limit} gives
\[
\frac{\Tr(F_k(f)^2)}{N_k}
\longrightarrow
\mu_{\rm eq}(f^2)-\mu_{\rm eq}(f)^2.
\]
Writing $F_k(f)=(u_{ij}^{(k)})$, the same bounds yield
\[
\begin{gathered}
0\leq\frac1{N_k}\sum_j(u_{jj}^{(k)})^2
\leq\frac{\Tr(F_k(f)^2)}{N_k},\qquad
\left|\frac{\Tr(F_k(f)F_k(f)^T)}{N_k}\right|
\leq\frac{\Tr(F_k(f)^2)}{N_k}.
\end{gathered}
\]
Thus the two remaining Bai--Yao profiles are bounded.
Every subsequence consequently admits a further subsequence
along which they converge. Theorem~\ref{T:Bai-Yao-mass-form}
then gives a centered Gaussian limit for $W_k$ along this
further subsequence. By Lemma~\ref{P:mass-covariance},
its variance is the corresponding subsequential limit
of $v_k(f)$.

If $\liminf_k v_k(f)>0$, this limiting variance is
strictly positive. Normalizing by $\sqrt{v_k(f)}$
and applying Slutsky's theorem, we obtain convergence
to $\mathcal N(0,1)$ along every such subsequence,
hence along the full sequence. If $v_k(f)\to v>0$,
the same argument gives $W_k\to\mathcal N(0,v)$
in distribution. For $v=0$, the conclusion follows
from $\E|W_k|^2=v_k(f)\to0$.

Since $\E r_k=1$ and
$\E(r_k-1)^2=(\nu_4-1)/N_k\to0$,
we have $r_k\to1$ in probability.
Identity~\eqref{e:ratio-identity}
and Slutsky's theorem therefore transfer these limits
to the self-normalized mass statistic centered at
$\mu_k^{\mathrm{Berg},\Delta}(f)$.
The convention on zero sections does not affect
the distributional limits, since
$\P(s_k=0)=\P(\xi=0)^{N_k}\to0$.

It remains to justify replacing
$\mu_k^{\mathrm{Berg},\Delta}(f)$ by
$\E\widehat{\cM}_{s_k}(f)$ in the centering.
Since the fluctuations are scaled by $\sqrt{N_k}$,
it suffices to show that the difference between
these two quantities is $o(N_k^{-1/2})$.
The convergence of $r_k$ in probability alone
does not provide this estimate.

We first estimate the difference between
the self-normalized mass statistic and
$\mu_k^{\mathrm{Berg},\Delta}(f)$.
To this end, let $G_k=\{r_k\geq1/2\}$. By
Lemma~\ref{L:mass-variance-comparison}
and the operator-norm bound on $F_k(f)$, we have
$\E|W_k|^2=v_k(f)\leq4C_\xi\|f\|_\infty^2$.
Moreover, Chebyshev's inequality gives
$\P(G_k^c)\leq4(\nu_4-1)/N_k$.

On $G_k$, the ratio identity~\eqref{e:ratio-identity}
and $r_k\geq1/2$ imply
$|\widehat{\cM}_{s_k}(f)
-\mu_k^{\mathrm{Berg},\Delta}(f)|
\leq2|W_k|/\sqrt{N_k}$.
On $G_k^c$, since both
$\widehat{\cM}_{s_k}$ and
$\mu_k^{\mathrm{Berg},\Delta}$ are probability
measures, the same difference is bounded by
$2\|f\|_\infty$. Consequently,
\[
\begin{aligned}
\E\bigl|\widehat{\cM}_{s_k}(f)
-\mu_k^{\mathrm{Berg},\Delta}(f)\bigr|^2
&\leq\frac4{N_k}\E|W_k|^2
+4\|f\|_\infty^2\P(G_k^c)=O(N_k^{-1}).
\end{aligned}
\]

To estimate the difference of expectations, we use
the relation
$r_k\widehat{\cM}_{s_k}(f)=\cM_{s_k}(f)$,
which also holds on $\{s_k=0\}$, where both sides
vanish. Since $\E r_k=1$ and
$\E\cM_{s_k}(f)=\mu_k^{\mathrm{Berg},\Delta}(f)$,
we obtain
$\E[r_k(\widehat{\cM}_{s_k}(f)
-\mu_k^{\mathrm{Berg},\Delta}(f))]=0$.
Hence, by Cauchy--Schwarz,
\[
\begin{aligned}
|\E\widehat{\cM}_{s_k}(f)
-\mu_k^{\mathrm{Berg},\Delta}(f)|
&=\left|\E\left[(1-r_k)
\bigl(\widehat{\cM}_{s_k}(f)
-\mu_k^{\mathrm{Berg},\Delta}(f)\bigr)\right]\right|\\
&\leq\bigl(\E(r_k-1)^2\bigr)^{1/2}
\bigl(\E|\widehat{\cM}_{s_k}(f)
-\mu_k^{\mathrm{Berg},\Delta}(f)|^2\bigr)^{1/2}\\
&=O(N_k^{-1}).
\end{aligned}
\]
We have therefore established
\begin{equation}\label{e:self-normalized-mean-comparison}
\E\widehat{\cM}_{s_k}(f)
-\E\cM_{s_k}(f)=O(N_k^{-1}).
\end{equation}
The implied constant depends only on $\|f\|_\infty$
and the coefficient law. In particular, this
difference is $o(N_k^{-1/2})$, so replacing
$\mu_k^{\mathrm{Berg},\Delta}(f)$ by
$\E\widehat{\cM}_{s_k}(f)$ does not affect either
distributional limit.

Finally, under the nondegeneracy
conditions~\eqref{e:CLT-nondegeneracy},
Lemma~\ref{L:mass-variance-comparison}
and~\eqref{e:centered-trace-limit} yield
\[
\liminf_{k\to\infty}v_k(f)
\geq c_\xi\bigl(
\mu_{\rm eq}(f^2)-\mu_{\rm eq}(f)^2
\bigr)>0
\]
whenever the equilibrium variance of $f$
is positive. 
\end{proof}

\begin{Remark}[Centering]\label{R:centering}
By~\eqref{e:self-normalized-mean-comparison},
the expectations of the dimension-normalized and
self-normalized mass statistics differ by
$O(N_k^{-1})$, so either may be used for centering
at the $\sqrt{N_k}$ scale. Centering instead at
$\mu_{\rm eq}(f)$ requires the additional condition
$\sqrt{N_k}(\mu_k^{\mathrm{Berg},\Delta}(f)
-\mu_{\rm eq}(f))\to0$, which does not follow
from weak convergence of the Bergman measures alone.
\end{Remark}

The covariance formula~\eqref{e:covariance-form}
simplifies when $\E\xi^2=0$ and $\E|\xi|^4=2$.
Under these conditions, $\beta=\kappa=0$, so only
the trace term remains. The following corollary
shows that both mass statistics then have the same
limiting Gaussian laws as for standard complex
Gaussian coefficients, independently of the
prescribed orthonormal bases. The coefficient
law itself need not be Gaussian.

\begin{Corollary}
\label{C:mass-fourth-matching}
Let $\xi$ be centered with $\E|\xi|^2=1$,
$\E\xi^2=0$, and $\E|\xi|^4=2$.
For every prescribed sequence of orthonormal bases
and every $f\in\mathcal C^0(K,\R)$,
\begin{equation}\label{e:fourth-matching-dimension}
\sqrt{N_k}\bigl(\cM_{s_k}(f)-\E\cM_{s_k}(f)\bigr)
\longrightarrow\mathcal N\bigl(0,\mu_{\rm eq}(f^2)\bigr),
\end{equation}
and
\begin{equation}\label{e:fourth-matching-self}
\sqrt{N_k}\bigl(\widehat{\cM}_{s_k}(f)
-\E\widehat{\cM}_{s_k}(f)\bigr)
\longrightarrow\mathcal N\bigl(0,
\mu_{\rm eq}(f^2)-\mu_{\rm eq}(f)^2\bigr).
\end{equation}
Both convergences are in distribution and include
the case of zero limiting variance.
\end{Corollary}

\begin{proof}
The moment assumptions give $\beta=\kappa=0$, hence
$\mathcal V_\xi(A,A)=\Tr(A^2)$ for every Hermitian $A$.
By~\eqref{e:CLT-Toeplitz-input} and
\eqref{e:centered-trace-limit},
\[
\begin{aligned}
N_k\Var\bigl(\cM_{s_k}(f)\bigr)
&=\frac{\Tr(A_k(f)^2)}{N_k}
\longrightarrow\mu_{\rm eq}(f^2),\\
v_k(f)
&=\frac{\Tr(F_k(f)^2)}{N_k}
\longrightarrow\mu_{\rm eq}(f^2)-\mu_{\rm eq}(f)^2.
\end{aligned}
\]
If $\mu_{\rm eq}(f^2)>0$,
Theorem~\ref{T:intro-exact-variance-CLT}
and Slutsky's theorem give
\eqref{e:fourth-matching-dimension}.
If $\mu_{\rm eq}(f^2)=0$, the centered statistic
on the left converges to zero in $L^2$
by the first variance limit.
The second variance limit and
Theorem~\ref{T:self-normalized-CLT}
give~\eqref{e:fourth-matching-self},
including the case of zero limiting variance.
\end{proof}

\begin{Remark}
Note that, the above corollary applies even to finitely supported
coefficient laws. For example, let $\xi=0$ with
probability $1/2$ and let $\xi$ take each of the values
$\sqrt2,-\sqrt2,i\sqrt2,-i\sqrt2$
with probability $1/8$.
This non-Gaussian law satisfies
$\E\xi=\E\xi^2=0$, $\E|\xi|^2=1$, and
$\E|\xi|^4=2$.
\end{Remark}

\subsection{Basis dependence on the projective line}

Although mass equidistribution is independent of
the prescribed orthonormal bases, the limiting
fluctuations may depend on their choice.
The following example illustrates this phenomenon
on $\mathbf P^1$ by comparing two ensembles with
the same coefficient law and test function,
but different orthonormal bases.

\begin{Example}\label{Ex:projective-bases}
Let $X=\mathbf P^1$, $L=\mathcal O(1)$,
and $h=h_{\rm FS}$ be the Fubini--Study metric.
Let $\mu=\omega_{\rm FS}$ be the normalized
Fubini--Study measure. Then $K=X$,
$N_k=k+1$, and $\mu_{\rm eq}=\mu$.
Consider the smooth real-valued function
\[
f(z)=\frac{|z|^2-1}{1+|z|^2},
\]
which extends across infinity. In the normalized monomial bases, for any coefficient
law satisfying~\eqref{e:CLT-moments},
\begin{equation}\label{e:projective-monomial-limit}
\sqrt{k+1}\,\cM_{s_k}(f)
\longrightarrow\mathcal N\bigl(0,(\nu_4-1)/3\bigr).
\end{equation}
For the same coefficient law, assuming in addition
that $\E\xi^2=0$, the orthonormal bases obtained
from the monomial bases by a normalized discrete
Fourier transform give
\begin{equation}\label{e:projective-Fourier-limit}
\sqrt{k+1}\,\cM_{s_k}(f)
\longrightarrow\mathcal N(0,1/3).
\end{equation}
Both convergences are in distribution. The corresponding
self-normalized statistics, centered at their
expectations, have the same limits. In particular,
for coefficients uniformly distributed on the unit
circle, the monomial statistic vanishes identically,
whereas the Fourier-transformed bases give
nondegenerate Gaussian fluctuations.
\end{Example}

\begin{proof}

In the affine chart $[1:z]$, let $e_0$ be
the standard frame of $\mathcal O(1)$.
With $dA$ denoting Euclidean area measure,
\[
|e_0|_{h_{\rm FS}}^2=\frac1{1+|z|^2},
\qquad
d\mu(z)=\frac{dA(z)}{\pi(1+|z|^2)^2}.
\]
The sections, expressed in this frame,
\[
S_{j,k}(z)
=\sqrt{(k+1)\binom{k}{j}}\,
z^j e_0^{\otimes k},
\qquad 0\leq j\leq k,
\]
form an orthonormal basis of
$H^0(\mathbf P^1,\mathcal O(k))$.

Distinct monomials are orthogonal by angular
integration, and the substitution $t=|z|^2$ gives
\[
\int_{\C}\frac{|z|^{2j}}{(1+|z|^2)^k}\,d\mu(z)
=\int_0^\infty\frac{t^j}{(1+t)^{k+2}}\,dt
=\frac{j!(k-j)!}{(k+1)!}.
\]
The binomial theorem gives $B_k(z)=k+1$, so the
normalized diagonal Bergman measure is exactly
$\mu$ for every $k$. The diagonal Bergman limit
therefore identifies $\mu_{\rm eq}$ with $\mu$.

Put $u=|z|^2/(1+|z|^2)$. After angular integration,
$\mu$ becomes Lebesgue measure $du$ on $[0,1]$,
and $f=2u-1$. In particular,
\[
\mu(f)=0,\qquad
\mu(f^2)=\int_0^1(2u-1)^2\,du=\frac13.
\]
Since $f$ is radial, $A_k(f)$ is diagonal in the
monomial basis. The beta-integral formula gives
its $j$th diagonal entry:
\[
\lambda_{j,k}
=2\frac{\int_0^1u^{j+1}(1-u)^{k-j}\,du}
        {\int_0^1u^j(1-u)^{k-j}\,du}-1
=\frac{2j-k}{k+2}.
\]
Thus $\Tr A_k(f)=0$ and
\begin{equation}\label{e:projective-trace-square}
\frac{\Tr(A_k(f)^2)}{k+1}
=\frac1{k+1}\sum_{j=0}^k\frac{(2j-k)^2}{(k+2)^2}
=\frac{k}{3(k+2)}\longrightarrow\frac13.
\end{equation}
For a diagonal Hermitian matrix $A$, the covariance
formula reduces to
$\mathcal V_\xi(A,A)=(\nu_4-1)\Tr(A^2)$.
If $\nu_4>1$, Theorem~\ref{T:intro-exact-variance-CLT}
and Slutsky's theorem, together with
\eqref{e:projective-trace-square}, give
\eqref{e:projective-monomial-limit}.
If $\nu_4=1$, then $|\xi|=1$ almost surely and
$\cM_{s_k}(f)=\Tr A_k(f)/(k+1)=0$
for every realization, so the degenerate
limit holds as well.

For the second choice of bases, let
$U_k=(u_{\ell j})_{0\leq\ell,j\leq k}$
be the normalized discrete Fourier matrix of
order $k+1$, defined by
\[
u_{\ell j}
=\frac1{\sqrt{k+1}}
e^{2\pi i\ell j/(k+1)}.
\]
Since $U_k$ is unitary, the sections
$\widetilde S_{j,k}
=\sum_{\ell=0}^k u_{\ell j}S_{\ell,k}$
form an orthonormal basis. The corresponding
Toeplitz matrix is
$\widetilde A_k=U_k^*A_k(f)U_k$.
As $|u_{\ell j}|^2=(k+1)^{-1}$,
all diagonal entries of $\widetilde A_k$
equal $\Tr A_k(f)/(k+1)=0$.
Unitary invariance also gives
$\Tr(\widetilde A_k^2)=\Tr(A_k(f)^2)$.
When $\beta=\E\xi^2=0$, the covariance formula
therefore yields
\[
\mathcal V_\xi(\widetilde A_k,\widetilde A_k)
=\Tr(\widetilde A_k^2)=\Tr(A_k(f)^2).
\]
Theorem~\ref{T:intro-exact-variance-CLT},
Slutsky's theorem, and
\eqref{e:projective-trace-square} now give
\eqref{e:projective-Fourier-limit}.
Here independent coefficients with the same law
are assigned in each basis, producing two distinct
random-section ensembles rather than two coordinate
representations of the same random section.

Finally, since $\Tr A_k(f)=0$ in both bases,
$\E\cM_{s_k}(f)
=\mu_k^{\mathrm{Berg},\Delta}(f)=0$.
Moreover, $r_k\to1$ in probability, and the
mean comparison~\eqref{e:self-normalized-mean-comparison}
gives
$\E\widehat{\cM}_{s_k}(f)=O((k+1)^{-1})$.
The ratio identity~\eqref{e:ratio-identity}
and Slutsky's theorem therefore give the same
limits for the self-normalized statistics centered
at their expectations.
\end{proof}

\subsection{Haar-typical orthonormal bases}

The preceding example shows that the limiting fluctuations
may depend on the choice of orthonormal basis. For
Haar-typical bases, however, a different picture emerges.
We show that, for almost every sequence of bases, the
two basis-dependent profiles appearing in the last
two terms of the covariance formula~\eqref{e:covariance-form}
have deterministic limits. This leads
to central limit theorems with explicit limiting variances:
the dimension-normalized variance retains a dependence
on the fourth moment of the coefficient law, whereas
the self-normalized variance is determined entirely
by the equilibrium measure.

Choose $k_0$ such that $N_k>0$ for $k\geq k_0$,
and fix a reference orthonormal basis $\{S_{j,k}\}_{j=1}^{N_k}$ of $\cH_k$ for each such $k$.
Every unitary matrix $V=(v_{\ell j})\in U(N_k)$
determines an orthonormal basis
\[
S_{j,k}^V=\sum_{\ell=1}^{N_k}v_{\ell j}S_{\ell,k},
\qquad 1\leq j\leq N_k.
\]
The corresponding Toeplitz matrix is denoted by
$A_k^V(f):=V^*A_k(f)V$.

Equip $\mathcal U:=\prod_{k\geq k_0}U(N_k)$
with the product of normalized Haar measures,
and denote the resulting probability measure by
$\P_{\rm H}$ and the corresponding expectation
by $\E_{\rm H}$. We write
$\mathbf U=(U_k)_{k\geq k_0}\in\mathcal U$
for a sequence of unitary matrices.
At each fixed degree $k$, normalized Haar
integration over $U(N_k)$ is denoted by $dV$.
Once a sequence of bases is fixed, the random
coefficients are independent copies of $\xi$
within each row, as in the fixed-basis model, and
no independence between rows is assumed.

\begin{Proposition}\label{T:Haar-basis-CLT}
For $\P_{\rm H}$-almost every
$\mathbf U=(U_k)_k\in\mathcal U$, the limits
\begin{equation}\label{e:Haar-profiles}
\begin{aligned}
\frac1{N_k}\sum_j\bigl(A_k^{U_k}(f)_{jj}\bigr)^2
&\longrightarrow\mu_{\rm eq}(f)^2,\qquad\frac1{N_k}\Tr\bigl(A_k^{U_k}(f)
(A_k^{U_k}(f))^T\bigr)
&\longrightarrow\mu_{\rm eq}(f)^2
\end{aligned}
\end{equation}
hold simultaneously for all $f\in\mathcal C^0(K,\R)$.
For each such fixed sequence of bases, every coefficient
law satisfying~\eqref{e:CLT-moments} gives
\begin{equation}\label{e:Haar-CLT}
\sqrt{N_k}\bigl(\cM_{s_k}(f)-\E\cM_{s_k}(f)\bigr)
\longrightarrow
\mathcal N\bigl(0,\mu_{\rm eq}(f^2)
+(\nu_4-2)\mu_{\rm eq}(f)^2\bigr)
\end{equation}
and
\begin{equation}\label{e:Haar-normalized-CLT}
\sqrt{N_k}\bigl(\widehat{\cM}_{s_k}(f)
-\E\widehat{\cM}_{s_k}(f)\bigr)
\longrightarrow
\mathcal N\bigl(0,\mu_{\rm eq}(f^2)
-\mu_{\rm eq}(f)^2\bigr)
\end{equation}
in distribution, for every $f\in\mathcal C^0(K,\R)$.
Both conclusions include zero limiting variance.
The full-measure set of basis sequences is independent
of the coefficient law.
\end{Proposition}

\begin{proof}
Fix $f\in\mathcal C^0(K,\R)$. For each
$k\geq k_0$, define the real-valued functions
$D_k,R_k$ on $U(N_k)$ by
\[
D_k(V):=\frac1{N_k}\sum_j
\bigl(A_k^V(f)_{jj}\bigr)^2,\qquad
R_k(V):=\frac1{N_k}
\Tr\bigl(A_k^V(f)(A_k^V(f))^T\bigr).
\]

Although $R_k(V)$ is real-valued, it need not
be nonnegative. For $\mathbf U=(U_k)_k\in\mathcal U$,
the limits in~\eqref{e:Haar-profiles} concern
$D_k(U_k)$ and $R_k(U_k)$ as $k\to\infty$.
We begin by computing their Haar expectations.

Fix $k$, and let $v_j$ denote the $j$th column
of $V\in U(N_k)$. Under normalized Haar measure,
$v_j$ is uniformly distributed on the unit sphere
in $\C^{N_k}$. The standard spherical fourth-moment
identity
\[
\int_{U(N_k)}
|(v_j)_a|^2|(v_j)_b|^2\,dV
=\frac{1+\delta_{ab}}{N_k(N_k+1)}
\]
(see~\cite[Chapter~2]{Meckes2019}),
together with unitary invariance, gives
\begin{equation}\label{e:Haar-profile-expectations}
\begin{aligned}
\E_{\rm H}[D_k(U_k)]
&=\int_{U(N_k)}D_k(V)\,dV=\frac{(\Tr A_k(f))^2+\Tr(A_k(f)^2)}
{N_k(N_k+1)}.
\end{aligned}
\end{equation}
Here we used
$A_k^V(f)_{jj}=v_j^*A_k(f)v_j$
and the fact that all columns have the same
distribution. The first equality follows
because the $k$th marginal of $\P_{\rm H}$
is normalized Haar measure on $U(N_k)$.

We next compute the expectation of $R_k$.
Since
\[
R_k(V)=D_k(V)+\frac1{N_k}
\sum_{i\ne j}\bigl(A_k^V(f)_{ij}\bigr)^2,
\]
it suffices to show that the off-diagonal
terms have zero Haar expectation.
For a diagonal unitary matrix
$Q=\operatorname{diag}(e^{i\theta_1},\ldots,
e^{i\theta_{N_k}})$, we have
$A_k^{VQ}(f)=Q^*A_k^V(f)Q$, and hence
\[
A_k^{VQ}(f)_{ij}
=e^{i(\theta_j-\theta_i)}A_k^V(f)_{ij}.
\]
When $i\ne j$, the square of this entry
can thus be multiplied by an arbitrary
unit complex phase. Right invariance of
Haar measure forces its integral to vanish.
Consequently,
$\E_{\rm H}[R_k(U_k)]=\E_{\rm H}[D_k(U_k)]$.

Finally, the Toeplitz trace limits give
$N_k^{-1}\Tr A_k(f)\to\mu_{\rm eq}(f)$
and
$N_k^{-1}\Tr(A_k(f)^2)\to\mu_{\rm eq}(f^2)$.
It follows from~\eqref{e:Haar-profile-expectations}
that the Haar expectations of both profiles
converge to $\mu_{\rm eq}(f)^2$.

We now pass from convergence of the Haar
expectations to almost-sure convergence.
Put $M=\|f\|_\infty$. For $V,W\in U(N_k)$, we have
\begin{equation}\label{e:Haar-HS-bound}
\|A_k^V(f)-A_k^W(f)\|_{\rm HS}
\leq2\|A_k(f)\|_{\rm op}\|V-W\|_{\rm HS}
\leq2M\|V-W\|_{\rm HS}.
\end{equation}
Since $\|A_k^V(f)\|_{\rm HS}\leq\sqrt{N_k}M$,
the Cauchy--Schwarz inequality gives
\[
\max\bigl\{|D_k(V)-D_k(W)|,\,
|R_k(V)-R_k(W)|\bigr\}
\leq\frac{4M^2}{\sqrt{N_k}}\|V-W\|_{\rm HS}.
\]
Here the estimate for $R_k$ follows from
$|\Tr(BC^T)|\leq\|B\|_{\rm HS}\|C\|_{\rm HS}$.
Thus both profiles are
$4M^2/\sqrt{N_k}$-Lipschitz with respect
to the Hilbert--Schmidt metric.

The concentration inequality for Haar measure on
$U(N_k)$, see~\cite[Theorem~5.17]{Meckes2019},
asserts that every real-valued $L$-Lipschitz
function $G$ on $U(N_k)$, with $L>0$, satisfies
\[
\P_{\rm H}\bigl(
|G(U_k)-\E_{\rm H}[G(U_k)]|>t
\bigr)
\leq 2\exp\left(
-\frac{(N_k-2)t^2}{24L^2}
\right),
\qquad t>0,
\]
for $N_k>2$.
Applying this inequality to $D_k$ and $R_k$,
whose Lipschitz constants are bounded by
$4M^2/\sqrt{N_k}$, we obtain, for every
$\varepsilon>0$,
\begin{equation}\label{e:Haar-concentration}
\begin{aligned}
&\P_{\rm H}\bigl(
|D_k(U_k)-\E_{\rm H}[D_k(U_k)]|>\varepsilon
\bigr)+\P_{\rm H}\bigl(
|R_k(U_k)-\E_{\rm H}[R_k(U_k)]|>\varepsilon
\bigr)
\leq4e^{-c_{f,\varepsilon}N_k^2},
\end{aligned}
\end{equation}
for all sufficiently large $k$, where
$c_{f,\varepsilon}>0$ depends only on
$\varepsilon$ and $\|f\|_\infty$ when $f\ne0$.
The case $f=0$ is immediate. Since $N_k\asymp k^n$, the right-hand side
is summable over $k$. The Borel--Cantelli
lemma, together with the convergence of
the Haar expectations, establishes
\eqref{e:Haar-profiles} for each fixed $f$.

To obtain a common full-measure set for all
$f\in\mathcal C^0(K,\R)$, apply the preceding
argument to the dense family $(f_\ell)_{\ell\geq1}$
used in~\eqref{e:weak-metric}, and intersect
the corresponding probability-one events.
The estimate
$\|A_k(f)-A_k(g)\|_{\rm op}
\leq\|f-g\|_\infty$
shows that the sum of the absolute changes
in the two profiles is bounded by
$2(\|f\|_\infty+\|g\|_\infty)\|f-g\|_\infty$,
uniformly in $k$ and $V$.
Since $\mu_{\rm eq}$ is a probability measure,
the same approximation controls
$|\mu_{\rm eq}(f)^2-\mu_{\rm eq}(g)^2|$.
Hence~\eqref{e:Haar-profiles} holds simultaneously
for every $f\in\mathcal C^0(K,\R)$ on a
single set of full $\P_{\rm H}$-measure.

Now fix a sequence $\mathbf U=(U_k)_k\in\mathcal U$
for which~\eqref{e:Haar-profiles} holds
simultaneously for all $f\in\mathcal C^0(K,\R)$.
The three normalized matrix profiles of $A_k^{U_k}(f)$
in Theorem~\ref{T:Bai-Yao-mass-form}
converge to $\mu_{\rm eq}(f^2)$,
$\mu_{\rm eq}(f)^2$, and
$\mu_{\rm eq}(f)^2$, respectively.
By the covariance formula~\eqref{e:covariance-form}, we have
\[
\begin{aligned}
\frac1{N_k}\mathcal V_\xi
\bigl(A_k^{U_k}(f),A_k^{U_k}(f)\bigr)
&\longrightarrow
\sigma_f^2:=\mu_{\rm eq}(f^2)
+(\nu_4-2)\mu_{\rm eq}(f)^2,
\end{aligned}
\]
where we used $|\beta|^2+\kappa=\nu_4-2$.
If $\sigma_f^2>0$, applying
Theorem~\ref{T:Bai-Yao-mass-form}
to $A_k^{U_k}(f)$ gives~\eqref{e:Haar-CLT}.
If $\sigma_f^2=0$, the exact covariance
formula gives convergence of the centered
quadratic forms to zero in $L^2$,
and the same conclusion follows.

For self-normalized mass, consider the trace-free matrices
\[
F_k^{U_k}(f):=A_k^{U_k}(f)
-\mu_k^{\mathrm{Berg},\Delta}(f)I_{N_k}.
\]
By~\eqref{e:centered-trace-limit}
and unitary invariance of the trace,
\[
\frac1{N_k}\Tr\bigl((F_k^{U_k}(f))^2\bigr)
\longrightarrow
\widehat{\sigma}_f^2
:=\mu_{\rm eq}(f^2)-\mu_{\rm eq}(f)^2.
\]
The other two normalized covariance profiles
vanish asymptotically. Indeed, since
$\mu_k^{\mathrm{Berg},\Delta}(f)
=N_k^{-1}\Tr A_k^{U_k}(f)$,
the diagonal-square and transpose-trace
profiles of $F_k^{U_k}(f)$ are exactly
$D_k(U_k)-\mu_k^{\mathrm{Berg},\Delta}(f)^2$
and
$R_k(U_k)-\mu_k^{\mathrm{Berg},\Delta}(f)^2$,
respectively. Both tend to zero by
\eqref{e:Haar-profiles} and
\eqref{e:diagonal-limit}.

Consequently, the covariance formula gives
\[
\frac1{N_k}\mathcal V_\xi
\bigl(F_k^{U_k}(f),F_k^{U_k}(f)\bigr)
\longrightarrow\widehat{\sigma}_f^2.
\]
Thus $v_k(f)\to\widehat{\sigma}_f^2$, and
Theorem~\ref{T:self-normalized-CLT}
gives~\eqref{e:Haar-normalized-CLT},
including the case $\widehat{\sigma}_f^2=0$.

The full-measure set was constructed solely
from the Toeplitz matrices and Haar measure,
and is therefore independent of the
coefficient law.
\end{proof}

The two limiting variances differ by the
contribution of the random total mass. Namely,
\[
\sigma_f^2-\widehat{\sigma}_f^2
=(\nu_4-1)\mu_{\rm eq}(f)^2.
\]
For Haar-typical bases, the self-normalized
limiting variance is thus independent of
the coefficient law, giving rise to a
universality phenomenon for self-normalized
fluctuations. The projective-line example
shows that this simplification need not
hold for specially chosen bases.

\appendix

\section{Complete convergence for quadratic forms}
\label{S:quadratic-moment}

Here we prove Lemma~\ref{P:complete-convergence}, which establishes
summability of the deviation probabilities for Hermitian
quadratic forms under the sharp moment condition
$\E|\xi|^{2+2/n}<\infty$.
The proof treats the diagonal and off-diagonal terms separately.
The diagonal terms are handled by a truncation argument, while
the off-diagonal terms are controlled by a second-moment estimate
when $n>1$ and a fourth-moment estimate when $n=1$.
The resulting bounds are uniform over matrices whose operator
norms are bounded independently of $k$.

\subsection{Summation and fourth-moment estimates}

We first record two summation estimates. If $N_k\asymp k^n$ and
$t\geq1$, then
\begin{align}
\sum_{k:N_k<t}N_k
&\leq C(1+t^{1+1/n}),
\label{e:dimension-tail-sum}\\
\sum_{k:N_k\geq t}N_k^{-a}
&\leq Ct^{1/n-a}, \qquad a>1/n.
\label{e:dimension-inverse-sum}
\end{align}
Indeed, for all sufficiently large $k$, we have
$ck^n\leq N_k\leq Ck^n$. If $N_k<t$, then
$k<(t/c)^{1/n}$, and the first estimate follows from
$\sum_{k\leq m}k^n\leq Cm^{n+1}$. Similarly, $N_k\geq t$
implies $k\geq(t/C)^{1/n}$, so comparison with the integral
of $x^{-na}$ gives
\[
\sum_{k:N_k\geq t}N_k^{-a}
\leq C\sum_{k\geq(t/C)^{1/n}}k^{-na}
\leq Ct^{1/n-a}.
\]
The finitely many initial indices are absorbed into the constants;
no monotonicity of $N_k$ is required.

Consequently, if $V$ is nonnegative and
$\E V^{1+1/n}<\infty$, Tonelli's theorem and
\eqref{e:dimension-tail-sum} give
\begin{equation}\label{e:tail-truncation-sum}
\begin{aligned}
\sum_kN_k\P(V>N_k)
&=\E\left[\sum_{k:N_k<V}N_k\right]\leq C(1+\E V^{1+1/n})<\infty.
\end{aligned}
\end{equation}
We will also use \eqref{e:dimension-inverse-sum} after interchanging
a sum over $k$ with an expectation involving truncated moments.

\begin{Lemma}[Fourth moment of an off-diagonal form]
\label{L:offdiagonal-fourth}
Suppose $\E\xi=0$, $\E|\xi|^2=1$ and $\E|\xi|^4<\infty$.
For every Hermitian matrix $A=(a_{ij})_{1\leq i,j\leq N}$ and
independent copies $\xi_1,\ldots,\xi_N$,
\begin{equation}\label{e:offdiag-fourth-endpoint}
\E\left|\sum_{i\ne j}a_{ij}\ov\xi_i\xi_j\right|^4
\leq C_\xi\bigl(\Tr(A^2)\bigr)^2,
\end{equation}
where $C_\xi$ depends only on the coefficient law.
\end{Lemma}

\begin{proof}
We first establish the fourth-moment bound for linear forms:
for arbitrary deterministic complex coefficients $c_j$,
\begin{equation}\label{e:linear-fourth-bound}
\E\left|\sum_jc_j\xi_j\right|^4
\leq C_\xi\left(\sum_j|c_j|^2\right)^2.
\end{equation}
To verify it, put $\beta=\E\xi^2$ and
$\kappa=\E|\xi|^4-2-|\beta|^2$. Expanding the fourth power,
independence and centering show that only terms in which every index
occurs at least twice contribute. The three pairings, with the
correction for the all-equal case, give
\[
\begin{aligned}
\E\left|\sum_jc_j\xi_j\right|^4
&=2\left(\sum_j|c_j|^2\right)^2
  +|\beta|^2\left|\sum_jc_j^2\right|^2+\kappa\sum_j|c_j|^4.
\end{aligned}
\]
The inequalities $|\sum_jc_j^2|\leq\sum_j|c_j|^2$ and
$\sum_j|c_j|^4\leq(\sum_j|c_j|^2)^2$ imply
\eqref{e:linear-fourth-bound}. The same bound holds with
$\xi_j$ replaced by $\ov\xi_j$.

Write $O=\sum_{i\ne j}a_{ij}\ov\xi_i\xi_j$.
Introduce independent Bernoulli variables
$\delta_1,\ldots,\delta_N$ of parameter $1/2$, independent of
the coefficients, and set $\Lambda=\{i:\delta_i=1\}$. Define
\[
B_\Lambda
=\sum_{i\in\Lambda,\ j\notin\Lambda}a_{ij}\ov\xi_i\xi_j.
\]
For each ordered pair $i\ne j$, the probability that
$i\in\Lambda$ and $j\notin\Lambda$ is $1/4$.
Thus $O=4\E_\delta B_\Lambda$ for fixed coefficient values,
and Jensen's inequality yields
\begin{equation}\label{e:Bernoulli-decoupling}
\E|O|^4\leq4^4\E_\delta\E|B_\Lambda|^4.
\end{equation}

The two index sets defining $B_\Lambda$ are disjoint, allowing us
to apply the linear-form bound successively to two independent
groups of coefficients. In particular, this argument requires only
the fourth moment of $\xi$, not its eighth moment. Fix $\Lambda$
and condition first on $(\xi_i)_{i\in\Lambda}$. For
$j\notin\Lambda$, put
$v_j=\sum_{i\in\Lambda}a_{ij}\ov\xi_i$, so that
$B_\Lambda=\sum_{j\notin\Lambda}v_j\xi_j$.
The variables $\xi_j$, $j\notin\Lambda$, are independent of the
coefficients $v_j$; hence \eqref{e:linear-fourth-bound} gives
\[
\E|B_\Lambda|^4
\leq C_\xi\E\left(\sum_{j\notin\Lambda}|v_j|^2\right)^2.
\]
Expanding the square and applying Cauchy--Schwarz to each pair,
\begin{align*}
\E\left(\sum_{j\notin\Lambda}|v_j|^2\right)^2
&=\sum_{j,\ell\notin\Lambda}\E(|v_j|^2|v_\ell|^2)\leq\left(\sum_{j\notin\Lambda}(\E|v_j|^4)^{1/2}\right)^2.
\end{align*}
A second application of \eqref{e:linear-fourth-bound}, this time
to $(\ov\xi_i)_{i\in\Lambda}$, yields
$(\E|v_j|^4)^{1/2}\leq C_\xi\sum_{i\in\Lambda}|a_{ij}|^2$.
Consequently,
\[
\E|B_\Lambda|^4
\leq C_\xi\left(
\sum_{i\in\Lambda,\ j\notin\Lambda}|a_{ij}|^2
\right)^2.
\]
Since $A$ is Hermitian, the sum of squared entries is at most
$\sum_{i,j}|a_{ij}|^2=\Tr(A^2)$. Substitution into
\eqref{e:Bernoulli-decoupling} proves the lemma.
\end{proof}

\subsection{Proof of the complete-convergence estimate}

\begin{proof}[Proof of Lemma~\ref{P:complete-convergence}]
Let $M=\sup_k\|A_k\|_{\rm op}<\infty$ and put $p=1+1/n$.
Decompose the centered, unnormalized quadratic form as
\[
z_k^*A_kz_k-\Tr A_k=D_k+O_k,
\]
where
\[
D_k=\sum_j a_{jj}^{(k)}Y_{j,k},\qquad
Y_{j,k}=|\xi_{j,k}|^2-1,\qquad
O_k=\sum_{i\ne j}a_{ij}^{(k)}\ov\xi_{i,k}\xi_{j,k}.
\]
The variables $Y_{j,k}$ are independent and centered
within each row, and all have the same distribution
as $Y=|\xi|^2-1$. The moment assumption implies
$\E|Y|^p<\infty$.

\smallskip\noindent
\emph{Truncation of the diagonal variables.}
For each row, define
\[
\begin{aligned}
\widetilde Y_{j,k}
&=Y_{j,k}\1_{\{|Y_{j,k}|\leq N_k\}}
-\E\bigl(Y\1_{\{|Y|\leq N_k\}}\bigr),\\
\widetilde D_k&=\sum_j a_{jj}^{(k)}\widetilde Y_{j,k}.
\end{aligned}
\]
Let $E_k=\{\max_{j\leq N_k}|Y_{j,k}|>N_k\}$.
By the union bound and \eqref{e:tail-truncation-sum}, applied
to $|Y|$,
\begin{equation}\label{e:truncation-summable}
\sum_k\P(E_k)
\leq\sum_kN_k\P(|Y|>N_k)<\infty.
\end{equation}

On $E_k^c$, no variable is removed by truncation. Thus the
only difference between $D_k$ and $\widetilde D_k$ is the
deterministic recentering term. Since
$|a_{jj}^{(k)}|\leq M$ and $\E Y=0$,
\begin{align*}
\frac{|D_k-\widetilde D_k|}{N_k}
&\leq M\left|\E\bigl(Y\1_{\{|Y|\leq N_k\}}\bigr)\right|\\
&=M\left|\E\bigl(Y\1_{\{|Y|>N_k\}}\bigr)\right|\\
&\leq M\E\bigl(|Y|\1_{\{|Y|>N_k\}}\bigr)
\longrightarrow0.
\end{align*}
The convergence follows from the integrability of $Y$ and
$N_k\to\infty$.

Fix $\varepsilon>0$. For all sufficiently large $k$, the preceding
bound is less than $\varepsilon/4$. On $E_k^c$, therefore,
$|D_k+O_k|>\varepsilon N_k$ implies
$|\widetilde D_k|>\varepsilon N_k/4$ or
$|O_k|>\varepsilon N_k/2$. It follows that
\begin{align}
\P(|D_k+O_k|>\varepsilon N_k)
&\leq\P(E_k)
+\P(|\widetilde D_k|>\varepsilon N_k/4)+\P(|O_k|>\varepsilon N_k/2).
\label{e:complete-probability-reduction}
\end{align}
By \eqref{e:truncation-summable}, it remains to show that the
last two probabilities are summable.

\smallskip\noindent
\emph{Case $n>1$.}
The truncated variables are centered and independent within each
row, so
\[
\E|\widetilde D_k|^2
=\sum_j(a_{jj}^{(k)})^2
\Var\bigl(Y\1_{\{|Y|\leq N_k\}}\bigr)
\leq M^2N_k\E\bigl(Y^2\1_{\{|Y|\leq N_k\}}\bigr).
\]
Chebyshev's inequality consequently gives
\[
\P(|\widetilde D_k|>\varepsilon N_k/4)
\leq\frac{16M^2}{\varepsilon^2N_k}
\E\bigl(Y^2\1_{\{|Y|\leq N_k\}}\bigr).
\]
To sum these bounds, Tonelli's theorem yields
\[
\begin{aligned}
\sum_k\frac1{N_k}\E\bigl(Y^2\1_{\{|Y|\leq N_k\}}\bigr)
&=\E\left[Y^2\sum_{k:N_k\geq|Y|}N_k^{-1}\right].
\end{aligned}
\]
Since $n>1$, the series $\sum_kN_k^{-1}$ converges, so the
contribution from $|Y|<1$ is finite. For $|Y|\geq1$,
\eqref{e:dimension-inverse-sum} with $a=1$ gives
$\sum_{k:N_k\geq|Y|}N_k^{-1}\leq C|Y|^{1/n-1}$.
It follows that
\[
\begin{aligned}
\sum_k\frac1{N_k}\E\bigl(Y^2\1_{\{|Y|\leq N_k\}}\bigr)
&\leq C+C\E\bigl(
|Y|^{1+1/n}\1_{\{|Y|\geq1\}}\bigr)<\infty.
\end{aligned}
\]
By Chebyshev's inequality, we conclude that
$\sum_k\P(|\widetilde D_k|>\varepsilon N_k/4)<\infty$.

For the off-diagonal term, the second-moment calculation in
\eqref{e:offdiag-L2}, before normalization by $N_k^2$, gives
$\E|O_k|^2\leq2\Tr(A_k^2)\leq2M^2N_k$.
Consequently,
\[
\P(|O_k|>\varepsilon N_k/2)
\leq\frac{8M^2}{\varepsilon^2N_k}.
\]
This is summable since $N_k\asymp k^n$ and $n>1$.
Thus both probabilities on the right of
\eqref{e:complete-probability-reduction} are summable in this case.

\smallskip\noindent
\emph{Case $n=1$.}
Here $p=2$, so $\E Y^2<\infty$ and $\E|\xi|^4<\infty$.
The series $\sum_kN_k^{-1}$ now diverges, so the preceding
second-moment estimates do not suffice. We use fourth moments
for $\widetilde D_k$ and $O_k$.

For independent centered real random variables $X_j$, expansion
of the fourth power gives
\[
\E\left(\sum_jX_j\right)^4
=\sum_j\E X_j^4+6\sum_{i<j}\E X_i^2\,\E X_j^2.
\]
Every other index pattern has an index occurring only once, and
its expectation vanishes by independence and centering. Since $A_k$
is Hermitian, its diagonal entries are real; thus the formula applies
to $X_j=a_{jj}^{(k)}\widetilde Y_{j,k}$.
For any real random variable $W$ with finite fourth moment,
$\E|W-\E W|^4\leq16\E|W|^4$, while
$\Var(W)\leq\E W^2$. Using $|a_{jj}^{(k)}|\leq M$, we obtain
\begin{equation}\label{e:endpoint-diagonal-fourth}
\E|\widetilde D_k|^4
\leq CM^4\left[
N_k\E\bigl(|Y|^4\1_{\{|Y|\leq N_k\}}\bigr)
+N_k^2(\E Y^2)^2\right].
\end{equation}
Although $\E|Y|^4$ may be infinite, each truncated fourth moment
is finite. Tonelli's theorem and \eqref{e:dimension-inverse-sum}
with $n=1$ and $a=3$ give
\[
\begin{aligned}
\sum_kN_k^{-3}\E\bigl(|Y|^4\1_{\{|Y|\leq N_k\}}\bigr)
&=\E\left[|Y|^4\sum_{k:N_k\geq|Y|}N_k^{-3}\right]\leq C+C\E\bigl(|Y|^2\1_{\{|Y|\geq1\}}\bigr)<\infty.
\end{aligned}
\]
Indeed, the contribution from $|Y|<1$ is finite because
$\sum_kN_k^{-3}<\infty$, while for $|Y|\geq1$ the estimate
\eqref{e:dimension-inverse-sum} gives
$\sum_{k:N_k\geq|Y|}N_k^{-3}\leq C|Y|^{-2}$.
Since also $\sum_kN_k^{-2}<\infty$, dividing
\eqref{e:endpoint-diagonal-fourth} by $N_k^4$ and summing
over $k$ gives
\[
\sum_k\frac{\E|\widetilde D_k|^4}{N_k^4}<\infty.
\]
The fourth-moment Markov inequality therefore yields
$\sum_k\P(|\widetilde D_k|>\varepsilon N_k/4)<\infty$.

Finally, Lemma~\ref{L:offdiagonal-fourth} and the operator-norm
bound give
\[
\E|O_k|^4
\leq C_\xi\bigl(\Tr(A_k^2)\bigr)^2
\leq C_\xi M^4N_k^2.
\]
Hence
\[
\P(|O_k|>\varepsilon N_k/2)
\leq\frac{C_{\xi,M,\varepsilon}}{N_k^2},
\]
which is summable because $N_k\asymp k$. This completes the case
$n=1$.

In both cases, the truncation errors are summable by
\eqref{e:truncation-summable}, and the remaining probabilities
are summable by the estimates above. Combining these facts with
\eqref{e:complete-probability-reduction} proves
\[
\sum_k\P\bigl(
|z_k^*A_kz_k-\Tr A_k|>\varepsilon N_k\bigr)<\infty.
\]
\end{proof}

\section{Mass equidistribution for nested coefficients}
\label{S:row-coupling}

The necessity of the higher moment condition in
Theorem~\ref{T:fixed-basis-dimension}\textup{(ii)}
depends on the coupling between coefficient rows.
For nested rows drawn from a single infinite iid
sequence, a finite second moment already suffices
for almost-sure mass equidistribution.

\begin{Proposition}[Nested coefficient rows]
\label{P:nested-rows}
In the geometric setting of Section~\ref{S:Geometry},
fix an orthonormal basis
$\{S_{1,k},\ldots,S_{N_k,k}\}$ for each $k$.
Let $Z_1,Z_2,\ldots$ be an iid sequence satisfying
\[
\E Z_1=0,\qquad
0<\sigma^2:=\E|Z_1|^2<\infty,
\]
and define
\[
s_k=\sum_{j=1}^{N_k}Z_jS_{j,k}.
\]
Then
\[
\cM_{s_k}\longrightarrow\sigma^2\mu_{\rm eq},
\qquad
\widehat{\cM}_{s_k}\longrightarrow\mu_{\rm eq}
\]
weakly almost surely. No moment above two is required.
For each prescribed sequence of bases, both limits hold
on a common event of probability one, which may depend
on the basis sequence.
\end{Proposition}

\begin{proof}
Rescaling the coefficients, we may assume $\sigma^2=1$.
Set $\xi_k=(Z_1,\ldots,Z_{N_k})^T$.

Fix $f\in\mathcal C^0(K,\R)$ and write $A_k=A_k(f)$.
By the quadratic-form identity~\eqref{e:mass-stat}
and the trace asymptotics~\eqref{e:diagonal-limit}, we know that
\[
\cM_{s_k}(f)=\frac1{N_k}\xi_k^*A_k\xi_k,
\qquad
\frac{\Tr A_k}{N_k}\longrightarrow\mu_{\rm eq}(f).
\]
Moreover, the operator-norm bound
\eqref{e:Toeplitz-op-bound} gives
$\|A_k\|_{\rm op}\leq\|f\|_\infty$.
Thus it remains to show that
\begin{equation}\label{e:nested-quadratic-convergence}
\frac{\xi_k^*A_k\xi_k-\Tr A_k}{N_k}
\longrightarrow0
\qquad\text{almost surely}.
\end{equation}

Write $M=\sup_k\|A_k\|_{\rm op}$.
For each positive integer $R$, truncate and recenter
the coefficients:
\[
Z_j^{(R)}
=Z_j\1_{\{|Z_j|\leq R\}}
-\E\bigl(Z_1\1_{\{|Z_1|\leq R\}}\bigr),
\qquad
\tau_R^2=\E|Z_1^{(R)}|^2.
\]
Let
$Z_k^{(R)}=(Z_1^{(R)},\ldots,Z_{N_k}^{(R)})^T$.
The truncated variables are bounded, centered, and
independent. If $\tau_R>0$, applying
Lemma~\ref{P:complete-convergence} to
$Z_j^{(R)}/\tau_R$ and using the Borel--Cantelli
lemma yields
\[
\frac{(Z_k^{(R)})^*A_kZ_k^{(R)}
-\tau_R^2\Tr A_k}{N_k}
\longrightarrow0
\qquad\text{almost surely}.
\]
For $\tau_R=0$, the conclusion is immediate.
Intersecting these probability-one events over
positive integers $R$, we obtain simultaneous
convergence for every $R$.

The common infinite sequence allows us to control
the truncation error by the strong law of large
numbers. Setting
$\rho_R^2=\E|Z_1-Z_1^{(R)}|^2$, we have,
on a further event of probability one,
simultaneously for every positive integer $R$,
\[
\frac{\|\xi_k\|^2}{N_k}\longrightarrow1,\qquad
\frac{\|Z_k^{(R)}\|^2}{N_k}\longrightarrow\tau_R^2,\qquad
\frac{\|\xi_k-Z_k^{(R)}\|^2}{N_k}
\longrightarrow\rho_R^2.
\]
Since $Z_1^{(R)}\to Z_1$ in $L^2$, we also have
$\tau_R\to1$ and $\rho_R\to0$ as $R\to\infty$.

For any Hermitian matrix $A$ and vectors $x,y$,
\[
|x^*Ax-y^*Ay|
\leq\|A\|_{\rm op}(\|x\|+\|y\|)\|x-y\|.
\]
Applying this inequality with $x=\xi_k$ and
$y=Z_k^{(R)}$, and using
$|\Tr A_k|/N_k\leq M$, we obtain
\[
\limsup_{k\to\infty}
\left|
\frac{\xi_k^*A_k\xi_k-\Tr A_k}{N_k}
\right|
\leq
M|\tau_R^2-1|+M(1+\tau_R)\rho_R.
\]
Letting $R\to\infty$ proves
\eqref{e:nested-quadratic-convergence}.

Together with the trace asymptotics, this gives
\[
\cM_{s_k}(f)\longrightarrow\mu_{\rm eq}(f)
\qquad\text{almost surely}.
\]
Applying this to the dense family $(f_\ell)_{\ell\geq1}$
used to define the metric~\eqref{e:weak-metric},
and intersecting the corresponding probability-one
events, we obtain
$\cM_{s_k}\to\mu_{\rm eq}$ weakly almost surely.

In particular, $r_k=\cM_{s_k}(K)\to1$ on this event,
so the sections are eventually nonzero. The identity
$\widehat{\cM}_{s_k}=r_k^{-1}\cM_{s_k}$ then gives
the self-normalized limit. Restoring the factor
$\sigma$ multiplies the dimension-normalized mass
by $\sigma^2$ and leaves the self-normalized mass
unchanged.
\end{proof}

To see the consequence for
Theorem~\ref{T:fixed-basis-dimension}\textup{(ii)},
consider a centered coefficient law satisfying
$\E|Z_1|^2=1$ but $\E|Z_1|^{2+2/n}=\infty$.
The proposition gives almost-sure mass equidistribution
for nested rows, whereas
Theorem~\ref{T:fixed-basis-dimension}\textup{(ii)}
rules it out for pairwise independent rows with the
same coefficient law. Thus the necessity of the higher
moment condition does not extend to arbitrary couplings
between rows.

\end{document}